\documentclass[11pt]{amsart}
\usepackage{geometry}                
\usepackage{graphicx}
\usepackage{amssymb}
\usepackage{epstopdf}
\usepackage{cancel}
\usepackage{appendix}
\newtheorem{theorem}{Theorem}[section]
\newtheorem{definition}{Definition}[section]

\newtheorem{remark}{Remark}[section]
\newtheorem{lemma}{Lemma}[theorem]
 \newtheorem{proposition}{Proposition}[section]
\newcommand{\F}{\widehat}
\newcommand{\n}{\|} 
\newcommand{\nn}{|\!|\!|}
\newcommand{\la}{\langle}
\newcommand{\ra}{\rangle}
\newcommand{\LL}{\mathcal L}
\newcommand{\R}{\mathbb{R}}

\title{Smoothing effect for time-degenerate Schr\"{o}dinger Operators}

\author{Serena Federico}
\thanks{This project has received funding from the European Union’s Horizon 2020 research and innovation programme under the Marie Sk\l odowska-Curie grant agreement No 838661 and No 777822}
\address{Department of Mathematics: Analysis, Logic and Discrete Mathematics,
Ghent University,
Krijgslaan 281, Building S8,
B 9000 Ghent,
Belgium}\email[Serena Federico ]{serena.federico@ugent.be}
\author{Gigliola Staffilani}
\thanks{G.S. was partially founded by the NSF grant No DMS-1764403}
\address{Department of Mathematics, Massachusetts Institute of Technology, 77 Massachusetts Avenue, Cambridge , MA, 02139, USA}\email[Gigliola Staffilani]{gigliola@math.mit.edu}

\begin{document}

\begin{abstract} 
In this work we consider an example of a linear time degenerate Schr\"odinger operator. We show that with the appropriate assumptions the operator satisfies a Kato smoothing effect. We also show that the solutions to the nonlinear initial value problems involving this operator and polynomial derivative nonlinearities are locally well-posed and their solutions also satisfy the same smoothing estimates as the linear solutions.
\end{abstract}

\maketitle
\section{Introduction}

In this paper we study the smoothing effect for time-degenerate Schr\"{o}dinger operators of the form
\begin{equation}
\label{P2}
\LL_\alpha=i\partial_t+t^\alpha\Delta_x+b(t,x)\cdot \nabla_x,
\end{equation}
where the coefficients $b(t,x)=(b_1(t,x),...,b_n(t,x))$ are \textit{complex valued} and satisfy suitable decay assumptions, and $\alpha>0$.

In the case $b(t,x)\equiv 0$, or when $b(t,x)=ct^\alpha$, where $c$ is a complex vector, we show below that standard Fourier analysis arguments can be applied, and the results will follow by application of more classical techniques.

In the more general situation  we deal with space-time variable coefficients, at least in the first order part, and we need to replace the standard use of the Fourier transform with the use of pseudo-differential calculus. This will allow us to obtain  smoothing estimates  for the linear operator \eqref{P2} and its non homogenous counterpart. Smoothing estimates for an operator such as \eqref{P2} where $\alpha=0$ are by now classical results, see for example  \cite{CR, PC, D1,KPV}.

Although the central part of this work is dedicated to  linear smoothing estimates, in the second part of the paper we also address   local well-posedness of the related Cauchy problem with derivative nonlinearities. As in the case of nondegenerate space variable coefficients Schr\"{o}dinger equations (see \cite{KPV1, KPRV, MMT} and references therein), also in our case local well-posedness relies heavily on the smoothing estimates proved in the first part of this work.

In proving smoothing estimates for the operator \eqref{P2}, the main problems  are given by the presence of the time degeneracy in the second order term and by the presence of the first order term $b(t,x)\cdot \nabla_x$. In particular, the time degeneracy has to be managed in order to apply a method similar to that of Mizohata \cite{M} and Doi \cite{D} to absorb the space-time variable coefficients  first order term through the application of the Sharp G\aa rding inequality. The presence of the term $b(t,x)\cdot \nabla_x$ affects the applicability of the Sharp G\aa rding inequality as well, and, in addition, determines, in general, a loss of derivatives. 

Due to these considerations  it is clear that conditions on $b(t,x)$ are necessary to control the behavior of the operator, and, specifically, conditions relating the coefficient $t^\alpha$ and the coefficients $b_j(t,x)$. It is also well known that even in the case $\alpha=0$ and $b(t,x)=b(x)$ conditions that are necessary (and sufficient)  need to be imposed  on $b(x)$ even for  local well-posedness to hold, see \cite{M}.

Our approach to overcome these problems in proving  the linear and non homogeneous smoothing effect estimates  is inspired by   techniques similar to those  in \cite{KPV1, KPRV}, which themselves are based on a method proposed by  Mizohata \cite{M} and Doi \cite{D}. The key point is the construction of a pseudo-differential operator $K$ possessing good (specific) properties with respect to the second order operator $t^\alpha\Delta_x$ and that permits to control the first order term $b(t,x)\cdot \nabla_x$. Since the leading coefficient depends only on the time-variable $t$, we can perform our analysis by keeping $t$ as a parameter provided that $b$ satisfies suitable conditions with respect to $t$ as well. This simplifies considerably the problem, since, this way, we can choose the operator $K$ related to $\Delta_x$ instead of $t^\alpha\Delta_x$.

The way the operator $K$ is introduced in the argument is by defining a new norm equivalent to the $H^s_x$-norm, and  in terms of $K$ itself.  Thanks to this new norm we will perform useful commutator estimates giving the smoothing effect with a gain of one derivative for the non homogeneous term with respect to the regularity of the initial data  (which is necessary to deal with derivative nonlinearities).

In order to extend the smoothing effect estimates we prove in the first part of this work  to the nonlinear problem with derivative nonlinearities we follow the work of Kenig, Ponce and Vega who   studied smoothing effect and  local well-posedness of the nonlinear Cauchy problem both for space-variable nondegenerate Schr\"{o}dinger operators and for space-variable nondegenerate ultrahyperbolic Schr\"{o}dinger operators (see  \cite{ KPV1, KPRV},  and also \cite{R, MMT} and references therein).

In the context of time degenerate Schr\"{o}dinger operators, at least in relation to Cauchy problems, to the best of our knowledge, the only result currently available is due to Cicognani and Reissig in \cite{CR}, in which they study the local well-posedness of the homogeneous Cauchy problem for time degenerate Schr\"{o}dinger operators of the same form as \eqref{P2} (or even more general) in Sobolev and Gevrey spaces. However, our approach is different from that used in \cite{CR}, in which neither the analysis of the smoothing effect nor the local well-posedness of the derivative nonlinear Cauchy problem is considered. 

In this paper we study the following nonlinear prototypes associated with \eqref{P2}:
\begin{equation}
\label{NLP1}
\left\{ \begin{array}{ll}
\LL_\alpha u= \pm u|u|^{2k}& \\
u(0,x)=u_0(x), & 
\end{array}\right .
\end{equation}
and 
\begin{equation}
\label{NLP2}
\left\{ \begin{array}{ll}
\LL_\alpha u= \pm  t^\beta \nabla u\cdot u^{2k},& \beta\geq\alpha>0,\\u(0,x)=u_0(x). & 
\end{array}\right .
\end{equation}

We remark immediately that although in this paper we explicitly prove results for \eqref{NLP2}, as mentioned again in Section 7, the same techniques we use here work when the nonlinearity is a  polynomial  $t^\beta P(u, \bar u, \nabla u, \nabla \bar u)$.
In the case \eqref{NLP1}, the proofs when\footnote{As mentioned above the case when $b(t,x)=ct^\alpha$ and $c$ a complex vector, can be treated like the case $b\equiv 0$.} $b\equiv 0$ and when $b\not\equiv 0$ are treated separately. This is done in order to show how classical technics can be applied to obtain the local well-posedness in the first case, while, in the more general case $b\not\equiv 0$, the use of pseudo-differential calculus will be needed to get the result.

In the IVP \eqref{NLP2}, we use a nonlinearity where the time factor appears. This choice is dictated by the fact that it allows for the application of  the weighted smoothing estimates at our disposal directly. In addition this type of nonlinearity does not give any limitation on the exponent $\alpha$, which will be any nonnegative real number. If $\beta>\alpha$ then a contraction mapping theorem based on the fact that in the analysis of the nonlinear term a power $t^{\beta-\alpha}$ appears makes the analysis pretty straightforward. 
On the other hand when $\beta=\alpha$ more care needs to be used. The result follows  from  the  combination of the technique used here for the IVP \eqref{NLP2} when $\beta>\alpha$, and that used in \cite{KPV1, KPRV} for the nondegenerate case (compare with the problem (1.3) in \cite{KPRV}) in order to remove the smallness of the initial data.  The main point consists in fact in modifying the functional space to which the solution belongs in a way that permits to obtain, via a mean value theorem, a time factor needed to apply the contraction argument. The norm to be used is the one used in \cite{KPV1, KPRV} where the time derivative of the solution is taken into account. Since the proof of this result follows the lines described above we shall omit the proof which is left to the reader.

We want to stress, once again, that the local solutions of the nonlinear problems \eqref{NLP1} and \eqref{NLP2}   above satisfy the weighted smoothing estimates stated below. 

Moreover, results concerning more general degenerate Schr\"{o}dinger operators can be obtained following the procedure described above. In fact, one might consider
$$\LL_\alpha=i\partial_t+\sum_{i,j=1}^nD_{x_i}a_{ij}(t,x)D_{x_j}+b(t,x)\cdot \nabla_x,$$
where $a_{ij}(t,x)$ is such that $a_{ij}(t,x)\sim t^\alpha a'_{i,j}(x)$, and, for all $i,j=1,...,n$, $a'_{i,j}\in C^\infty_b(\mathbb{R}^n)$, $\{a'_{i,j}(x)\}_{i,j=1,...,n}$ is real valued and positive definite, and  $|\partial^\gamma_x a_{i,j}(t,x)|\lesssim t^\alpha \la x\ra^{-|\gamma|-\sigma+1}$ for any multiindex $\gamma\in \mathbb{N}^n$, where $\sigma$ is dictated by the conditions satisfied by $b(t,x)$ (see \eqref{condb} below). By combining our technique and that in \cite{KPRV} one should be able to prove the weighted smoothing effect, and, as a consequence, that the associated IVP (both linear and nonlinear similar to the ones above) is locally well-posed.
\subsection*{Notations}
We recall here, briefly, some notations used throughout the paper.

We use the notation $A\lesssim B$ to indicate that there exists an absolute constant $c>0$ such that  $A\leq cB$. 
We shall denote by $\Lambda^s$ the pseudo-differential operator of order $s$ whose symbol is given by $\Lambda^s(\xi)=\la \xi\ra^s=(1+|\xi|^2)^{s/2}$.

Below we often use mixed norm spaces. For example  $L^p_xL^q_t([0,T]\times\R^n), \, 1\leq p,q\leq \infty$ is the space of functions $f(x,t)$ that are in $L^q$ in time on the interval $[0,T]$ and are $L^p$ in space. The norm is taken in the right to left order. In a similar manner we define the spaces $L^p([0,T];H^s(\mathbb{R}^n)), \, 1\leq p\leq \infty$ of functions that are $L^p$ in time and in the Sobolev space $H^s(\mathbb{R}^n)$ in space. 
As remarked in \cite{KPV1}, also in this case if the  vector $b(t,x)=b(x)$ is a real smooth enough function,  then standard energy method gives well-posedness. On the other hand in this paper we are not just concerned  with 
 well-posedness, but  more importantly with the proof of smoothing estimates. 

We shall now summarize the main results of the paper. 
\subsection*{Case $b\equiv0$} 
We state here the results when $b\equiv0$. We will remark later that the same results hold when $b=ct^\alpha$, where $c$ is an imaginary vector.
\subsection*{Smoothing effect estimates}
Let $W_\alpha(t,s)$ be the operator defined as in \eqref{Wts}.
\begin{theorem} \label{teo0}
Let $W_\alpha(t):=W_\alpha(t,0)$, with $\alpha> 0$, then\\

If $n=1$ for all $f\in L^2(\mathbb{R})$,

\begin{equation}
\label{sh1}
\sup_x \n t^{\alpha/2}D_x^{1/2}W_\alpha(t)f\n^2_{L^2_t([0,T])}\lesssim  \n f\n^2_{L^2(\mathbb{R})};
\end{equation}
\\

If $n\geq2$, on denoting by $\{Q_\beta\}_{\beta\in\mathbb{Z}^n}$ the family of non overlapping cubes of unit size such that $\mathbb{R}^n=\bigcup_{\beta\in\mathbb{Z}^n}Q_\beta$, then for all $f\in L^2_x(\mathbb{R}^n)$,

\begin{equation}
\label{sh2}
\sup_{\beta\in\mathbb{Z}^n}\left( \int_{Q_\beta} \int_{0}^{T} | t^{\alpha/2}D_{x}^{1/2}W_\alpha(t)f(x)|^2 dt\,dx\right)^{1/2}\lesssim  \n f\n_{L^2(\mathbb{R}^n)},
\end{equation}
where $D_x^\gamma f(x)= (|\xi|^\gamma \widehat{f}(\xi))^\vee(x).$ 
\end{theorem}

\begin{theorem} \label{teo01}
Let  $n=1$ and $g\in L^1_xL^2_t([0,T]\times\mathbb{R})$, then 
\begin{equation}
\label{ish0}
\n D_x^{1/2}\int_{\mathbb{R_+}}t^{\alpha/2}W_\alpha(0,t)g(t)dt\n_{L^2_x(\mathbb{R})}\lesssim \n g\n_{L^1_{x}L^2_{t}(\mathbb{R}\times [0,T])},
\end{equation}
and, for all $g\in L^1_{t}L^2_{x}([0,T]\times \mathbb{R})$, 
\begin{equation}
\label{ish}
\n t^{\alpha/2}D_x^{1/2}\int_0^tW_\alpha(t,\tau)g(\tau)d\tau\n_{L^\infty_x(\mathbb{R})L^2_t{([0,T])}}\lesssim  \n g\n_{L^1_{t}L^2_{x}([0,T]\times \mathbb{R})}.
\end{equation}
If $n\geq 2$, denoting by $\{Q_\beta\}_{\beta\in\mathbb{Z}^n}$ a family of non overlapping cubes of unit size such that $\mathbb{R}^n=\bigcup_{\beta\in\mathbb{Z}^n}Q_\beta$, then,  for all $g\in L^1_{t}L^2_{x}([0,T]\times \mathbb{R}^n)$,

\begin{equation}
\label{ishn}
\sup_{\beta\in\mathbb{Z}^n}\left( \int_{Q_\beta}\left \n  t^{\alpha/2}D_{x}^{1/2} \int_{0}^{t}  W_\alpha(t,\tau)g(\tau) d\tau\right\n^2_{L^2_t([0,T])}\,dx\right)^{1/2}\lesssim  \n g\n_{L^1_{t}L^2_{x}([0,T]\times \mathbb{R}^n)},
\end{equation}
\end{theorem}

\subsection*{Local well-posedness}
\begin{theorem}\label{teo001}
Let $k\geq 1$, then \eqref{NLP1} is locally well-posed in $H^s$ for $s>n/2$ and its solution satisfies smoothing estimates.
\end{theorem}

\subsection*{Case $b\not\equiv0$} 
Let us consider the IVP
\begin{equation}
\label{IVP2}
\left\{ \begin{array}{ll}
\partial_tu=it^\alpha\Delta_xu +ib(t,x)\cdot \nabla_xu+f(t,x)& \\
u(0,x)=u_0(x). & 
\end{array}\right .
\end{equation}

The results true in this case are true in general, that is, also when $b\equiv 0$.
\subsection*{Smoothing effect estimates}
\begin{theorem} \label{thb}
Let $u_0\in H^s(\mathbb{R}^n), s\in \mathbb{R}$. 
Assume that, for all $j=1,...,n$, $b_j$ is such that $b_j\in C([0,T],C_b^\infty(\mathbb{R}^n))$ and there exists $\sigma>1$ such that
\begin{equation}
\label{condb}
| \mathsf{Im}\,\partial_x^\gamma b_j(t,x)|, | \mathsf{Re}\,\partial_x^\gamma b_j(t,x)|\lesssim t^\alpha \langle x \rangle^{-\sigma-|\gamma|},\quad x\in\mathbb{R}^n,
\end{equation}
and denote by $\lambda(|x|):=\langle x \rangle^{-\sigma}$.\\
Then
\begin{itemize}
  \item[(i)] If  $f\in L^1([0,T];H^s(\mathbb{R}^n))$ then the IVP \eqref{IVP2} has a unique solution $u\in C([0,T];H^s(\mathbb{R}^n))$
  and there exist positive constants $C_1, C_2$ such that
  \begin{equation*}
\sup_{0\leq t\leq T} \n u(t)\n_s\leq C_1 e^{C_2(\frac{T^{\alpha+1}}{\alpha+1}+T)}\left( \n u_0\n_s+\int_0^T \n f(t)\n_s dt\right);
\end{equation*}\\
  
  \item[(ii)] If  $f\in (L^2[0,T];H^s(\mathbb{R}^n))$ then the IVP \eqref{IVP2} has a unique solution $u\in C([0,T];H^s(\mathbb{R}^n))$
  and there  exist two positive constants $C_1, C_2$ such that
 \begin{equation*}
\sup_{0\leq t\leq T} \n u(t)\n^2_s +\int_0^T \int_{\mathbb{R}^n}  t^\alpha \left| \Lambda^{s+1/2} u \right|^2 \lambda(|x|) dx\,dt
\leq C_1 e^{C_2(\frac{T^{\alpha+1}}{\alpha+1}+T)}\left( \n u_0\n^2_s+\int_0^T \n f(t)\n^2_s dt\right);
\end{equation*}\\

  \item[(iii)] If  $\Lambda^{s-1/2} f\in L^2([0,T]\times \mathbb{R}^n;t^{-\alpha}\lambda(|x|)^{-1}dtdx)$ then the IVP \eqref{IVP2} has a unique solution $u\in C([0,T];H^s(\mathbb{R}^n))$
  and there exist positive constants $C_1, C_2$ such that
 \begin{equation*}
 \sup_{0\leq t\leq T} \n u(t)\n^2_s +\int_0^T \int_{\mathbb{R}^n}  t^\alpha \left| \Lambda^{s+1/2} u \right|^2 \lambda(|x|) dx\,dt
 \end{equation*}
$$\leq C_1 e^{C_2\frac{T^{\alpha+1}}{\alpha+1}}\left( \n u_0\n^2_s+\int_0^T  \int_{\mathbb{R}^n} t^{-\alpha} \lambda(|x|)^{-1}   \left| \Lambda^{s-1/2} f \right|^2 dx\, dt\right).$$

\end{itemize}
\end{theorem}
Above we abbreviated the norm $\|f\|_{H^s(\R^n)}=:\|f\|_s$.

\subsection*{Local well-posedness}
\begin{theorem}\label{teo1}
Let $\LL_\alpha$ be such that condition \eqref{condb} is satisfied. Then the IVP \eqref{NLP1} is locally well posed in $H^s$ for $s>n/2$.
\end{theorem}

\begin{theorem}\label{teo2}
Let $\LL_\alpha$ be such that condition \eqref{condb} is satisfied with $\sigma=2N$ (thus $\lambda(|x|)=\la x \ra^{-2N}$) for some $N\geq 1$,
and $s>n+4N+3$ such that $s-1/2\in 2\mathbb{N}$.  Let $H^s_\lambda:=\{u_0\in H^s(\mathbb{R}^n); \lambda(|x|)u_0\in H^s(\mathbb{R}^n)\}$, then the IVP \eqref{NLP2} with $\beta\geq\alpha>0$, is locally well posed in $H^s_\lambda$.
\end{theorem}

We remark that in all our local well-posedness results we do not focus on the optimal index $s$ for the $H^s$ regularity of the initial data. Instead our goal in these theorems is to show that the nonlinear solutions enjoy the same type of smoothing estimates as the linear ones.

We conclude this introduction by giving the plan of the paper.
In Section 2 we start the analysis of the case $b\equiv 0$, and we remark how $b=ct^\alpha$, with   $c$ a constant imaginary vector, can be analyzed  in a similar manner. First we show that the solution (which is explicit) to the linear  initial value problem involving $\LL_\alpha$ can be written via  the use of a two-parameter family of unitary operators (similarly to the standard case). Afterwards, by using standard Fourier analysis methods we derive homogeneous time-weighted smoothing estimates (where the weight depends on the degeneracy) by reducing the case $\alpha\neq 0$ (degenerate case) to the nondegenerate case $\alpha=0$. This is the first part of Theorem \ref{teo0}.

In Section 3 we show that the Duhamel's principle still applies in our context and prove inhomogeneous time-weighted smoothing estimates with a gain of 1/2 derivative with respect to the initial data. This is Theorem \ref{teo01}.

In Section 4 we state the local well-posedness result for the nonlinear Cauchy problem in the case $b\equiv 0$ in which we consider a nonlinearity of the form $u|u|^{2k}$ with $k\geq 1$, and we prove that smoothing estimates propagate to these nonlinear solutions.  This is Theorem \ref{teo001}.

In Section 5 we start the analysis of the more general case $b\not\equiv 0$ (but the results are true when $b\equiv0$ as well). We state the hypothesis on the coefficients $b_j$ and prove the local smoothing effect estimates in weighted Sobolev spaces. Additionally, we prove the local well-posedness of the linear Cauchy problem by means of the smoothing estimates. This is Theorem \ref{thb}.

In Section 6 we analyze the local well-posedness of the nonlinear Cauchy problem in presence of two different nonlinearities. First we consider nonlinearities of the form $u|u|^{2k}$, $k\geq 1$. This is Theorem \ref{teo1}.   Afterwards, we consider derivative nonlinearities of the form $t^\beta\nabla u \cdot u^{2k}$, where $\nabla u:=\mathsf{div}(u)$. This is Theorem \ref{teo2}.

Finally Section 7 contains some final remarks. 

\section{The case $\LL_\alpha=i\partial_t+t^\alpha\Delta_x$. Homogeneous Smoothing properties}
We start with the analysis of the homogeneous Cauchy  problem

\begin{equation}
\label{HCP}
\left\{ \begin{array}{ll}
\partial_tu=it^\alpha\Delta_xu & \\
u(s,x)=u_s(x), & 
\end{array}\right .
\end{equation}
where $0\leq s<t\leq T$, $x\in\mathbb{R}^n$ and $u_0$ is at least in $L^2(\mathbb{R}^n)$.

Observe that, by application of the Fourier transform with respect to the space variable, we get 

\begin{displaymath}
\label{ }
\left\{ \begin{array}{ll}
\partial_t \F{u}(t,\xi)=-it^\alpha |\xi|^2 \F{u}(t,\xi)&\\
\F{u}(s,\xi)=\F{u}_s(\xi), & 
\end{array}\right .
\end{displaymath}
 whose solution at time $t\leq T$ is given by
 $$\F{u}(t,\xi)=e^{-i\frac{t^{\alpha+1}-s^{\alpha+1}}{\alpha+1}|\xi|^2}\F{u}_s(\xi),$$
and finally, by Fourier inversion formula,
\begin{equation}
\label{uh}
u(t,x)=\int_{\mathbb{R}^n}e^{-i(\frac{t^{\alpha+1}-s^{\alpha+1}}{\alpha+1}|\xi|^2-x\cdot\xi)}\F{u}_s(\xi)d\xi.
\end{equation}
Formula \eqref{uh}, giving the solution of the homogeneous problem at $0<t\leq T$ starting at time $s<t$, can be written as
\begin{equation}
\label{uh2}
u(t,x)= W_\alpha(t,s)u_s(x):=e^{i\frac{t^{\alpha+1}-s^{\alpha+1}}{\alpha+1}\Delta_x}u_s(x):=\frac{(\alpha +1)^{n/2}}{(i(t^{\alpha+1}-s^{\alpha+1}))^{n/2}}e^{i(\alpha+1)\frac{|\cdot|^2}{t^{\alpha+1}-s^{\alpha+1}}}\ast u_s(x).
\end{equation}
Therefore throughout the paper we shall use the notation
\begin{equation}
\label{Wts}
W_\alpha(t,s):=e^{i\frac{t^{\alpha+1}-s^{\alpha+1}}{\alpha+1}\Delta_x}, \quad \forall s,t\in[0,T],
\end{equation}
to indicate the operator defined as in \eqref{uh2}.
Note that $\{W_\alpha(t,s)\}_{s,t\in[0,T]}$ is a two-parameter family of unitary operators and that, for any given $(t,s)$, $W_\alpha(t,s)$ is the ''solution operator'' of the IVP  \eqref{HCP}, that is, $u(t,x)= W_\alpha(t,s)u_s(x)$ is the solution at time $t$ of \eqref{HCP} starting at time $s$.
Moreover the following properties hold:
\begin{description}
  \item[(i)] $W_\alpha (t,t)=I;$
  \item[(ii)] $W_\alpha(t,s)= W_\alpha(t,r)W_\alpha(r,s)$ for every $s,t,r\in[0,T]$;
   \item[(iii)] $W_\alpha(t,s)\Delta_x u=\Delta_xW_\alpha(t,s) u$.
\end{description}

In particular we are interested in solutions starting at time $s=0$ (where the operator is degenerate), that, by \eqref{Wts}, will be given by $u(t,x)=W_\alpha(t,0)u_0(x)$.
In  what follows we shall often use the notation $W_\alpha(t,0)=:W_\alpha(t)$.

Note that, if $\alpha=0$ and $t>0$,
\begin{equation}
\label{W0}
W_0(t)=e^{i t \Delta_x},
\end{equation}
which is the standard Schr\"odinger semigroup.

By \eqref{Wts} we can easily see the first property of $W_\alpha$, namely
\begin{equation}
\label{ }
\n W_\alpha(t,s)u_s\n_{H^s_x}=\n u_s\n_{H^s_x}.
\end{equation}
The second property of the operator $W_\alpha$ is given by the  local smoothing result of Theorem \ref{teo0} that we prove below.

\proof[Proof of Theorem \ref{teo0}]
First note that \eqref{sh1} and \eqref{sh2} are true when $\alpha=0$ in \eqref{HCP}, that is, when $W_\alpha(t)=W_0(t)=e^{it\Delta_x}$ is the standard Schr\"{o}dinger semigroup (see, for instance, \cite{KPV}).
Then it suffices to prove that 
\begin{equation}
\label{proof1}
\n t^{\alpha/2}D_x^{1/2}W_\alpha(t)f\n^2_{L^2_t([0,T])}=C_\alpha\n D_x^{1/2}W_0(t)f\n^2_{L^2_t([0,T'])},
\end{equation}
where $C_\alpha$ is a positive constant depending on $\alpha$ and $T'>0$, since then the result will follow directly from the standard case $\alpha=0$.
We then reduce the proof to the proof of \eqref{proof1}.

We have that
\begin{equation}
\label{proof2}
\n t^{\alpha/2}D_x^{1/2}W_\alpha(t)f\n^2_{L^2_t([0,T])}=\int_0^T\left| t^{\alpha/2}\int_{\mathbb{R}^n}e^{-i\left(t^{\alpha+1}|\xi|^2/(\alpha+1)-x\cdot\xi\right)}|\xi|^{1/2}\F{f}(\xi)d\xi  \right|^2dt,
\end{equation}
then we apply the change of variables $t^{\alpha+1}/(\alpha+1)=s$, $t=c_\alpha s^{1/(\alpha+1)},\,c_\alpha=(\alpha+1)^{1/(\alpha+1)}$, $dt=c'_\alpha s^{-\alpha/(\alpha+1)}ds$  with $c'_\alpha=c_\alpha/(\alpha+1)$, and get
$$\eqref{proof2}=\int_0^{\frac{T^{\alpha+1}}{\alpha+1}}c_\alpha^\alpha \cancel{s^{\alpha/(\alpha+1)}}\left|\int_{\mathbb{R}^n}e^{-i\left(s|\xi|^2-x\cdot\xi\right)}|\xi|^{1/2}\F{f}(\xi)d\xi  \right|^2c'_\alpha \cancel{s^{-\alpha/(\alpha+1)}}ds$$
$$=c_\alpha^\alpha c'_\alpha \int_0^{\frac{T^{\alpha+1}}{\alpha+1}} \left|\int_{\mathbb{R}^n}e^{-i\left(t|\xi|^2-x\cdot\xi\right)}|\xi|^{1/2}\F{f}(\xi)d\xi  \right|^2dt$$
$$= C_\alpha \int_0^{\frac{T^{\alpha+1}}{\alpha+1}}|D_x^{1/2}W_0(t)f(x)|^2dt$$
$$=C_\alpha \n D_x^{1/2}W_0(t)f\n^2_{L_t^2([0,T^{\alpha+1}/(\alpha+1)])},$$
which gives the result.
Finally, by application of the smoothing estimates for $W_0(t)=e^{it\Delta_x}$, we get \eqref{sh1} and \eqref{sh2} (see \cite{KPV}, Corollary 2.2).
\endproof

\begin{remark}
Now we consider the case $b(x,t)=ct^\alpha$, where $c$ is imaginary. We note that in his case
\begin{equation}
\label{Wtsc}
W_\alpha(t,s):=
e^{i\frac{t^{\alpha+1}-s^{\alpha+1}}{\alpha+1}(\Delta_x +c\cdot \nabla_x)} s,t\in[0,T],
\end{equation}
and 
$$\widehat{W_\alpha(t,0)u_0}(\xi)=e^{-i\frac{t^{\alpha+1}}{\alpha+1}(|\xi|^2-ic\cdot \xi)},$$
and the argument proceeds exactly as above.
\end{remark}

\section{The case $-i\partial_t+t^\alpha\Delta_x$. Inhomogeneous smoothing properties}
We consider here the inhomogeneous Cauchy problem

\begin{equation}
\label{ICP}
\left\{ \begin{array}{ll}
\partial_tu=it^\alpha\Delta_xu +f(t,x)& \\
u(0,x)=u_0(x), & 
\end{array}\right .
\end{equation}
where $(t,x)\in [0,T]\times\mathbb{R}^n$, $u_0$ is at least in $L^2(\mathbb{R}^n)$ and $f$ is at least in $L^2([0,T]\times\mathbb{R}^n)$, and we prove that Duhamel's formula still applies in this case.

\begin{proposition}
Let $u_0\in L^2_x(\mathbb{R}^n)$, then the solution at time $t>0$ of the IVP (initial value problem) \eqref{ICP} is given by
\begin{equation}
\label{D}
u(t,x)=W_\alpha(t)u_0(x)+\int_0^t W_\alpha(t,\tau)f(\tau,x)d\tau.
\end{equation}
\end{proposition}

\proof
To prove \eqref{D} we will prove that the solution of the IVP \eqref{ICP} is given by
\begin{equation}
\label{D1}
u(t,x)= u_1(t,x)+\int_0^t u_2(\tau,x)d\tau,
\end{equation}
where $u_1$ and $u_2$ are the solutions of 

\begin{eqnarray*}
\text{IVP1}= \left\{ \begin{array}{ll}
\partial_tu_1=it^\alpha\Delta_xu & \\
u_1(0,x)=u_0(x), & 
\end{array}\right .& &\text{IVP2}=\left\{ \begin{array}{ll}
\partial_tu_2=it^\alpha\Delta_xu & \\
u_2(\tau,x)=f(\tau,x), & 
\end{array}\right .
\end{eqnarray*}
respectively.

By \eqref{uh} we have that, for $t$ fixed and $\tau<t$, 
$$u_1(t,x)=W_\alpha(t)u_0(x),$$
and
$$u_2(t,x)=W_\alpha(t,\tau)f(\tau,x):=e^{i(\frac{t^{\alpha+1}-{\tau}^{\alpha+1}}{\alpha+1})\Delta_x}f(\tau,x).$$\\
We now suppose that $u(t,x)$ is the solution at time $t$ of \eqref{ICP}, and show that $u$ is exactly given by \eqref{D1}. 
Note that, by using the fact that $W_\alpha(t,s)$ commutes with $\Delta_x$, then if $u$ solves \eqref{ICP} we get
$$\frac{d}{dt}W_\alpha(0,t)u(t,x)=\frac{d}{dt}\left( e^{-i \frac{t^{\alpha+1}}{\alpha+1}\Delta_x}u(t,x)\right)$$
$$=-it^\alpha \Delta_x W_\alpha(0,t)u(t,x)+W_\alpha(0,t)\frac{d}{dt}u(t,x)$$
$$=-it^\alpha \Delta_x W_\alpha(0,t)u(t,x)+W_\alpha(0,t)\left(it^\alpha\Delta_xu(t,x)+f(t,x) \right)$$
$$=\cancel{-it^\alpha \Delta_x W_\alpha(0,t)u(t,x)}+\cancel{it^\alpha W_\alpha(0,t)\Delta_x u(t,x)}+ W_\alpha(0,t)f(t,x),$$
which gives
$$\frac{d}{dt}W_\alpha(0,t)u(t,x)=W_\alpha(0,t)f(t,x).$$
We then integrate the last equality 
$$\int_0^t\frac{d}{ds}W_\alpha(0,s)u(s,x)ds=\int_0^tW_\alpha(0,\tau)f(\tau,x)d\tau$$
and find
$$W_\alpha(0,t)u(t,x)=W_\alpha(0,0)u_0(x)+\int_0^tW_\alpha(0,\tau)f(\tau,x)d\tau,$$
which gives, by applying $W_\alpha(t):=W_\alpha(t,0)$ on both sides, and recalling that $W_\alpha(t,t)=I$,
$$u(t,x)=W_\alpha(t)u_0(x)+W_\alpha(t)\int_0^tW_\alpha(0,\tau)f(\tau,x)d\tau$$
$$=W_\alpha(t)u_0(x)+\int_0^tW_\alpha(t,\tau)f(\tau,x)d\tau$$
$$=u_1(t,x)+\int_0^t u_2(\tau,x)d\tau.$$
Therefore, if $u$ solves \eqref{ICP} it is of the form \eqref{D1}, that is, $u$ is given by formula \eqref{D}.

We conclude the proof by verifying that $u$ given by \eqref{D} satisfies \eqref{ICP}.

We have
$$\frac{d}{dt} u(t,x)\underset{\eqref{D}}{=}it^\alpha \Delta_x W_\alpha(t)u_0(x)+\int_0^t it^\alpha\Delta_x W_\alpha(t,\tau)f(\tau,x)\, d\tau+f(t,x)$$
$$=it^\alpha \Delta_x\left( W_\alpha(t)u_0(x)+\int_0^t W_\alpha(t,\tau)f(\tau,x)\, dt\right)+f(t,x)$$
$$=it^\alpha \Delta_x u(t,x)+f(t,x),$$
thus $u$ satisfies \eqref{ICP}.
\endproof

We are now ready to prove Theorem \ref{teo01}.
\proof[Proof of Theorem \ref{teo01}]

The first inequality \eqref{ish0} follows directly from \eqref{sh1} by duality.

As regards \eqref{ish}, denoting $L^p_x:=L^p_x(\mathbb{R}^n)$, we have
$$ \n t^{\alpha/2}D_x^{1/2}\int_0^t W_\alpha(t,\tau)g(\tau)d\tau\n_{L^\infty_x L^2_t([0,T])}$$
$$\leq\left\n \left(\int_0^T\left| \int_0^T \left|t^{\alpha/2}D_x^{1/2} W_\alpha(t,\tau)g(\tau)\right|\,d\tau\right|^2dt\right)^{1/2}\right\n_{L_x^\infty}$$
$$\underset{\text{Minkowski}}{\leq} \left\n \int_0^T\left(  \int_0^T  \left | t^{\alpha/2}D_x^{1/2}W_\alpha(t,\tau)g(\tau)\right|^2 dt\right)^{1/2}d\tau\right\n_{L^\infty_x}$$
$$=\int_0^T \left\n t^{\alpha/2}D_x^{1/2}W_\alpha(t,0) (W_\alpha(0,\tau)g(\tau))\right\n_{L^\infty_x L^2_t([0,T])}d\tau$$
$$\underset{\text{by \eqref{sh1}}}{\leq}\int_0^T \n W_\alpha(0,\tau)g(\tau)\n_{L_x^2}d\tau= \n g\n_{L^1_{t}([0,T])L_x^2},$$
which gives \eqref{ish}.

As for \eqref{ishn} we first observe that
$$\left\n  t^{\alpha/2}D_{x}^{1/2} \int_{0}^{t}  W_\alpha(t,\tau)g(\tau) d\tau\right\n_{L^2_t([0,T])}$$
$$\leq \left\n \int_{0}^{T}  \left|t^{\alpha/2}D_{x}^{1/2} W_\alpha(t,\tau)g(\tau)\right| d\tau\right\n_{L^2_t([0,T])}$$
$$\underset{\text{Minkowski}}{\leq } \int_0^T \n t^{\alpha/2}D_{x}^{1/2} W_\alpha(t,0)W_\alpha(0,\tau)g(\tau)\n_{L^2_t([0,T])}d\tau,$$
therefore

$$\left(\int_{Q_\beta}\left \n  t^{\alpha/2}D_{x}^{1/2} \int_{0}^{t}  W_\alpha(t,\tau)g(\tau) d\tau\right\n^2_{L^2_t([0,T])}\,dx\right)^{1/2}$$
$$\leq \left[  \int_{Q_\beta} \left(\int_0^T \n t^{\alpha/2}D_{x}^{1/2} W_\alpha(t)W_\alpha(0,\tau)g(\tau)\n_{L^2_t([0,T])}d\tau\right)^2 dx\right]^{1/2}$$
$$\underset{\text{Minkowski}}{\leq}  \int_0^T\left ( \int_{Q_\beta} \n t^{\alpha/2}D_{x}^{1/2} W_\alpha(t,0)W_\alpha(0,\tau)g(\tau)\n^2_{L^2_t([0,T)]} dx\right)^{1/2}d\tau.$$
We then apply the  $\sup_{\beta\in\mathbb{Z}^n}$ on both the RHS and the LHS of the latter inequality and get

$$\sup_{\beta\in\mathbb{Z}^n}\left(\int_{Q_\beta}\left \n  t^{\alpha/2}D_{x}^{1/2} \int_{0}^{t}  W_\alpha(t,\tau)g(\tau) d\tau\right\n^2_{L^2_t([0,T])}\,dx\right)^{1/2}$$
$$\leq \sup_{\beta\in\mathbb{Z}^n} \left( \int_0^T\left ( \int_{Q_\beta} \n t^{\alpha/2}D_{x}^{1/2} W_\alpha(t,0)W_\alpha(0,\tau)g(\tau)\n^2_{L^2_t([0,T])} dx\right)^{1/2}d\tau\right)$$
$$ \leq \int_0^T \sup_{\beta\in\mathbb{Z}^n} \left( \int_{Q_\beta} \n t^{\alpha/2}D_{x}^{1/2} W_\alpha(t,0)W_\alpha(0,\tau)g(\tau)\n^2_{L^2_t([0,T])} dx\right)^{1/2}d\tau$$

$$\underset{\text{by \eqref{sh2}}}{\leq}  \int_0^T \n W_\alpha(0,\tau)g(\tau)\n_{L^2_x(\mathbb{R}^n)} d\tau$$
$$=\int_0^T \n g(\tau)\n_{L^2_x(\mathbb{R}^n)}d\tau= \n g\n_{L^1_{t}([0,T])L^2_x(\mathbb{R}^n)}$$
which concludes the proof.

\endproof
As remarked for the homogeneous problem, also in the inhomogeneous one we considered above we can take $b(x,t)=t^\alpha c$ where $c$ is an imaginary vector and obtain the same results. 

\section{Local well-posedness for the nonlinear Cauchy problem}

Let us now consider the nonlinear initial value problem
\begin{equation}
\label{NCP}
\left\{ \begin{array}{ll}
\partial_tu=it^\alpha\Delta_xu \pm u|u|^{2k}& \\
u(0,x)=u_0(x), & 
\end{array}\right .
\end{equation}
where $k\geq 1$ is a positive  integer. We shall prove that \eqref{NCP} is locally well posed in $H^s(\mathbb{R}^n)$, for $s>n/2$.\\

 In the sequel we shall use the notation $W_\alpha(t,0)=:W_\alpha(t)$.

\begin{definition} We say that the IVP (\ref{NCP}) is locally well-posed (l.w.p) in $H^s(\mathbb{R}^n)$ if for any ball $B$ in the space $H^s(\mathbb{R}^n)$ there exist a time $T$ and a Banach space of functions $X\subset L^\infty([0,T],H^s(\mathbb{R}^n))$ such that for each initial datum $u_0 \in B$ there exists a unique solution $u \in X \subset C ([0,T], H^s (\mathbb{R}^n ))$ for the integral equation

$$u(x, t) = W_\alpha(t)u_0 + \int_ 0^t W_\alpha(t ,\tau)|u|^{2k}u(\tau) d\tau.$$
Furthermore the map $u_0 \mapsto u$ is continuous as a map from $H^s(\mathbb{R}^n)$ into $C([0,T],H^s(\mathbb{R}^n))$.
\end{definition}

To prove the local well-posedness we shall use the following result.

\begin{lemma}
Let $g(u)=u |u|^{2k}$ and $s,l$ positive integers with $l\leq s$ and $s>n/2$. Then
\begin{equation}
\label{e1}
\n g(u)\n_{H^s}\lesssim \n u\n_{H^s}^{2k+1},
\end{equation}
\begin{equation}
\label{e2}
\n g(u)-g(v)\n_{L^2}\lesssim (\n u\n_{H^s}^{2k}+\n v\n_{H^s}^{2k})\n u-v\n_{L^2},
\end{equation}
\begin{equation}
\label{e3}
\n g^{(l)}(u)-g^{(l)}(v)\n_{L^\infty}\lesssim (\n u\n_{H^s}^{2k-l}+\n v\n_{H^s}^{2k-l})\n u-v\n_{H^s},
\end{equation}
\begin{equation}
\label{e4}
\n g(u)-g(v)\n_{H^s}\lesssim (\n u\n_{H^s}^{2k}+\n v\n_{H^s}^{2k})\n u-v\n_{H^s}.
\end{equation}
\end{lemma}

We are now ready to prove Theorem \ref{teo001}.

\proof[Proof of Theorem \ref{teo001}]
 The proof is based on the standard fixed point argument. For convenience we shall assume that the nonlinear term is given by $+u|u|^{2k}$ but the proof applies with no modification in the focusing case.
We assume first that $n=1$.
Let $X$ be the following metric space
$$X=\{ u:[0,T]\times\mathbb{R}\rightarrow \mathbb{C}; \n t^{\alpha/2}D^{1/2+s}_x u\n_{L^\infty_xL^2_t([0,T])}<\infty,\, \n u\n_{L^\infty_t([0,T])H^s_x}<\infty \},$$
equipped with the distance
$$d(u,v)=\n t^{\alpha/2}D^{1/2+s}_x (u-v)\n_{L^\infty_xL^2_t([0,T])}+\n u-v\n_{L^\infty_t([0,T])\dot{H}^s_x}+\n u-v\n_{L^\infty_t([0,T])L^2_x}$$
in which $\dot{H}^s_x$ stands for the homogeneous Sobolev space, and consider
$$\Phi:X\rightarrow X,\quad \Phi(u)=W_\alpha(t)u_0+\int_0^tW_\alpha(t,\tau) u|u|^{2k}(\tau)d\tau.$$
We now prove that $\Phi$ is a contraction, since then the result follows by application of the fixed point theorem.

We have that
$$\n \Phi(u)\n_X\leq \n W_\alpha(t)u_0\n_X+\n \int_0^tW_\alpha(t,\tau) u|u|^{2k}(\tau)d\tau\n_X,$$
and consider the two terms on the RHS separately.\\

For the homogeneous term we have
$$\n W_\alpha(t)u_0\n_X=\n t^{\alpha/2}D^{1/2+s}_x W_\alpha(t) u_0\n_{L^\infty_xL^2_t([0,T])}+\n W_\alpha(t)u_0\n_{L^\infty_t([0,T])\dot{H}^s_x}+ \n W_\alpha(t) u_0\n_{L^\infty_t([0,T])L^2_x}$$
$$\underset{\text{by\eqref{sh1}}}{\leq} \n D^s_xu_0\n_{L^2_x}+ 2\n u_0\n_{H^s_x}$$
$$\leq 3 \n u_0\n_{H^s}.$$

For the nonlinear term we get
$$\n \int_0^tW_\alpha(t,\tau) u|u|^{2k}(\tau)d\tau\n_X=  \n t^{\alpha/2}D^{1/2+s}_x \int_0^tW_\alpha(t,\tau) u|u|^{2k}(\tau)d\tau \n_{L^\infty_xL^2_{[0,T]}}$$
$$+\n \int_0^tW_\alpha(t,\tau)u|u|^{2k}(\tau)d\tau\n_{L^\infty_t([0,T])\dot{H}^s_x}+\n \int_0^tW_\alpha(t,\tau) u|u|^{2k}(\tau)d\tau\n_{L^\infty_t([0,T])L^2_x}$$
$$\underset{\eqref{ish}}{\leq} \n D^s_xu|u|^{2k}\n_{L^1_{t}([0,T])L^2_x}+\n \int_0^tW_\alpha(t,\tau) D^s_x u|u|^{2k}(\tau)d\tau\n_{L^\infty_t([0,T])L^2_x}$$
$$+\n \int_0^tW_\alpha(t,\tau) u|u|^{2k}(\tau)d\tau\n_{L^\infty_t([0,T])L^2_x}.$$
Note that
$$\n D^s_x u|u|^{2k}\n_{L^1_{t}([0,T])L^2_x}$$
$$\leq T \n u|u|^{2k}\n_{L^\infty_{t}([0,T])H^s_x}$$
$$ \underset{\text{by \eqref{e1}}}{\leq} C T \n u\n^{2k+1}_{L^\infty_{t}([0,T])H^s_x}$$
$$\leq C T \n u\n^{2k+1}_X,$$
and that
$$\n \int_0^tW_\alpha(t,\tau) D^s_x u|u|^{2k}(\tau)d\tau\n_{L^\infty_t([0,T])L^2_x}=\n W_\alpha(t) \int_0^tW_\alpha(0,\tau) D^s_x u|u|^{2k}(\tau)d\tau\n_{L^\infty_t([0,T])L^2_x}$$
$$\leq \left(\int_\mathbb{R}\left| \int_0^T \left|W_\alpha(0,\tau) D^s_x u|u|^{2k}(\tau)\right|d\tau\right|^2 dx\right)^{1/2}$$
$$\leq \int_0^T  \left(\int_\mathbb{R}\left|W_\alpha(0,\tau) D^s_x u|u|^{2k}(\tau)\right|^2 dx\right)^{1/2} d\tau$$
$$= \n D^s_x u|u|^{2k}\n_{L^1_t([0,T])L^2_x}$$
$$\leq T \n u|u|^{2k}\n_{L^\infty_{t}([0,T])H^s_x}$$
$$ \leq C T \n u\n^{2k+1}_{L^\infty_{t}([0,T])H^s_x}$$
$$\leq C T \n u\n^{2k+1}_X.$$
Similarly
$$\n \int_0^tW_\alpha(t,\tau) u|u|^{2k}(\tau)d\tau\n_{L^\infty_t([0,T])L^2_x}$$
$$\leq \n D^s_x u|u|^{2k}\n_{L^1_{t}([0,T])L^2_x}$$
$$\leq C T \n u\n^{2k+1}_X,$$
therefore, by application of the previous inequalities we get
$$\n \int_0^tW_\alpha(t,\tau) D^s_x u|u|^{2k}(\tau)d\tau\n_X \leq  C T \n u\n_X^{2k+1}.$$
Putting together the previous estimates we have
$$\n \Phi(u)\n_X\leq 3 \n u_0\n_{H^s_x}+ C_1T \n u\n_X^{2k+1}$$
which gives that $\Phi:X\rightarrow X$. 

Let $R= 6 \n u_0\n_{H^s_x}$ and consider $B_R=\{ u\in X;\, \n u\n_X \leq R\}\subset X$. Then, by choosing $T$ such that $ C_1 T R^{2k}\leq 1/2$, we get that, for all $u\in B_R$
$$ \n \Phi(u)\n_X\leq R/2 +  C_1T R^{2k+1}\leq R,$$
which gives that $\Phi$ sends $B_R$ in $B_R$.

We now prove that $\Phi$ is a contraction.

Let $u,v\in B_R$, and consider
$$ \n \Phi(u)-\Phi(v)\n_X=\n \int_0^t W_\alpha(t,\tau)(u|u|^{2k}(\tau)-v|v|^{2k}(\tau))d\tau\n_X$$
$$= \n t^{\alpha/2}D_x^{s+1/2} \int_0^t W_\alpha(t,\tau)\left(u|u|^{2k}(\tau)-v|v|^{2k}(\tau)\right)d\tau\n_{L^\infty_xL^2_t([0,T])}$$
$$+\n \int_0^t W_\alpha(t,\tau)\left(u|u|^{2k}(\tau)-v|v|^{2k}(\tau)\right)d\tau\n_{L^\infty_t([0,T])\dot{H}^s_x}$$
$$+\n \int_0^t W_\alpha(t,\tau)\left(u|u|^{2k}(\tau)-v|v|^{2k}(\tau)\right)d\tau\n_{L^\infty_t([0,T])L^2_x}$$
$$= I+II+III.$$
By using the same estimates used before we have that 
$$I\leq \n D^s (u|u|^{2k}-v|v|^{2k}) \n_{L^1_{t}([0,T])L^2_x}$$
$$\leq T \n u|u|^{2k}-v|v|^{2k} \n_{L^\infty_{t}([0,T])H^s_x}$$
$$\underset{\text {by \eqref{e4}}}{\leq} CT (\n u\n_{L^\infty_{t}([0,T])H^s_x}^{2k}+\n v\n_{L^\infty_{t}([0,T])H^s_x}^{2k})\n u-v\n_{L^\infty_{t}([0,T])H^s_x}$$
$$\leq CT (\n u\n_X^{2k}+\n v\n_X^{2k})\n u-v\n_X$$
$$\leq C TR^{2k}\n u-v\n_X.$$

For $II$ we estimate the $L^\infty_t([0,T])\dot{H}^s_x$ - norm as before and get, once more,
$$II\leq \n D^s (u|u|^{2k}-v|v|^{2k}) \n_{L^1_{t}([0,T])L^2_x}$$
$$\leq C TR^{2k}\n u-v\n_X,$$
while, once again,
$$ III\leq \n u|u|^{2k}-v|v|^{2k}\n_{L^1_{t}([0,T])L^2_x}$$
$$\leq C TR^{2k}\n u-v\n_X.$$
We then have
$$ \n \Phi(u)-\Phi(v)\n_X\leq C_2 TR^{2k}\n u-v\n_X,$$
where $T$ was chosen in such a way that $ C_1 T R^{2k}\leq 1/2$. We then choose the time $T$ such that $T=\min\{\frac{1}{C_1 R^{2k}}, \frac{1}{C_2 R^{2k}}\}$, and conclude that $\Phi$ is a contraction. The result then follows by fixed point arguments.

We now assume that $n>1$. We proceed like above and  we  consider the space

$$X=\{ u:[0,T]\times\mathbb{R}^n\rightarrow \mathbb{C}; \nn t^{\alpha/2}D^{s+1/2}_x u\nn_T<\infty,\, \n u\n_{L^\infty_{[0,T]}H^s_x}<\infty \},$$
where 
$$ \nn \cdot\nn_T=\sup_{\beta\in\mathbb{Z}^n} \n \cdot\n_{L^2_x(Q_\beta)L^2_t([0,T])},$$
and 
$$d_X(u,v)= \nn t^{\alpha/2}D^{s1/2}_x (u-v)\nn_T+ \n u-v\n_{L^\infty_t([0,T])\dot{H}^s_x}+ \n u-v\n_{L_t^\infty([0,T]) L^2_x}.$$

Once again we consider 
$$\Phi:X\rightarrow X,\quad \Phi(u)=W_\alpha(t)u_0+\int_0^tW_\alpha(t,\tau) u|u|^{2k}(\tau)d\tau,$$
and prove that $\Phi$ is a contraction.

We have that
$$\n \Phi(u)\n_X\leq \n W_\alpha(t)u_0\n_X+\n \int_0^tW_\alpha(t,\tau) u|u|^{2k}(\tau)d\tau\n_X,$$
and we estimate the two terms on the RHS separately.

For the homogeneous term we have 
$$\n W_\alpha(t)u_0\n_X= \nn t^{\alpha/2}D^{s+1/2}_x W_\alpha(t)u_0\nn_T+\n W_\alpha(t)u_0\n_{L^\infty_t([0,T])\dot{H}^s_x}+\n W_\alpha(t)u_0\n_{L^\infty_t([0,T])L^2_x}$$
$$\leq  \nn t^{\alpha/2}D^{1/2}_x W_\alpha(t)D^s_x u_0\nn_T+2\n W_\alpha(t)u_0\n_{L^\infty_t([0,T])H^s_x}$$
$$= \sup_{\beta\in\mathbb{Z}^n} \n  t^{\alpha/2}D^{1/2}_x W_\alpha(t)D^s_x u_0\n_{L^2(Q_\beta)L^2_t([0,T])}+2\n W_\alpha(t)u_0\n_{L^\infty_t([0,T])H^s_x}$$
$$\underset{\eqref{ishn}}{\leq} \n D^s_x u_0\n_{L^2(\mathbb{R}^n)}+2\n u_0\n_{H^s_x}$$
$$\leq 3\n u_0\n_{H^s_x}.$$

For the inhomogeneous term we have 
\begin{equation}
\label{interm}
\n \int_0^tW_\alpha(t,\tau) u|u|^{2k}(\tau)d\tau\n_X\leq  \underbrace{\nn \int_0^tW_\alpha(t,\tau) u|u|^{2k}(\tau)d\tau\nn_T}_{(\ref{interm}.1)}
\end{equation}
$$+ \underbrace{ 2\n \int_0^tW_\alpha(t,\tau) u|u|^{2k}(\tau)d\tau\n_ {L^\infty_t([0,T])H^s_x}}_{(\ref{interm}.2)},$$
where
$$(\ref{interm}.1)= \sup_{\beta\in\mathbb{Z}^n} \n  t^{\alpha/2}D^{1/2}_x \int_0^t W_\alpha(t,\tau)D^s_x u|u|^{2k}(\tau)d\tau\n_{L^2(Q_\beta)L^2_t([0,T])}$$
$$\underset{\eqref{ishn}}{\leq} \n D^s_x u|u|^{2k}\n_{L^1_{t}([0,T])L^2(\mathbb{R}^n)}\leq T \n u|u|^{2k}\n_{L^\infty_{t}([0,T])H^s(\mathbb{R}^n)}$$
$$\leq T \n u\n^{2k+1}_{L^\infty_{t}([0,T])H^s(\mathbb{R}^n)}\leq C T \n u\n_X^{2k+1},$$
and
$$(\ref{interm}.2)=\sup_{t\in[0,T]} \left\n \int_0^tW_\alpha(t,\tau) u|u|^{2k}(\tau)d\tau\ \right\n_{H^s(\mathbb{R}^n)}$$
$$=\sup_{t\in[0,T]} \left\n \int_0^tW_\alpha(t,\tau) (1+D_x^{2s})^{1/2}u|u|^{2k}(\tau)d\tau\ \right\n_{L^2_x}$$
$$\leq \sup_{t\in[0,T]} \left(\int_{\mathbb{R}^n}\left| \int_0^t |{W}_\alpha(t,\tau)  (1+D_x^{2s})^{1/2} u|u|^{2k}(\tau)|d\tau\right|^2\ dx \right)^{1/2}$$
$$\underset{\text{Minkowski}}{\leq} \sup_{t\in[0,T]}  \left( \int_0^t \n   (1+D_x^{2s})^{1/2}u|u|^{2k}(\tau)\n_{L^2(\mathbb{R}^n)}d\tau\right)$$
$$= \int_0^T \n  u|u|^{2k}(\tau)\n_{H^s(\mathbb{R}^n)}d\tau$$
$$=\n  u|u|^{2k}\n_{L^1_{t}([0,T])H^s(\mathbb{R}^n)}\leq T \n  u|u|^{2k}\n_{L^\infty_{t}([0,T])H^s(\mathbb{R}^n)}$$
$$ \leq T \n u\n^{2k+1}_{L^\infty_{t}([0,T])H^s(\mathbb{R}^n)}\leq C T \n u\n_X^{2k+1}.$$

Therefore, putting together the previous estimates, we get
$$\n \Phi(u)\n_X\leq 3\n u_0\n_{H^s}+C_1T\n u\n_X^{2k+1},$$
hence $\Phi:X\rightarrow X$.

Let  now $R= 6\n u_0\n_{H^s(\mathbb{R}^n)}$ and $T$ such that $C_1T R^{2k}\leq 1/2$, then  $\n \Phi(u)\n_X\leq R$ for all $u\in B_R$, that is, $\Phi$ sends $B_R$ into $B_R$.

We now end the proof by showing that $\Phi$ is a contraction.

Let $u,v\in B_R$, then
$$ \n \Phi(u)-\Phi(v)\n_X=\n \int_0^t W_\alpha(t,\tau)(u|u|^{2k}(\tau)-v|v|^{2k}(\tau))d\tau\n_X$$
$$\leq \nn t^{\alpha/2}D_x^{s+1/2}\int_0^tW_\alpha(t,\tau)(u|u|^{2k}(\tau)-v|v|^{2k}(\tau))d\tau\nn_T$$
$$+ \n \int_0^t W_\alpha(t,\tau)(u|u|^{2k}(\tau)-v|v|^{2k}(\tau))d\tau\n_{L^\infty_t([0,T])H^s_x(\mathbb{R}^n)}$$
$$=I+II,$$
where, by \eqref{e4} and the procedure used before,

$$I, II\leq \n u|u|^{2k}-v|v|^{2k}\n_{L^1_{t}([0,T])H^s_x(\mathbb{R}^n)}$$
$$\leq CT (\n u\n^{2k}_{L^\infty_{t}([0,T])H^s_x(\mathbb{R}^n)}+\n v\n^{2k}_{L^\infty_{t}([0,T])H^s_x(\mathbb{R}^n)}) \n u-v\n^{2}_{L^\infty_{t}([0,T])H^s_x(\mathbb{R}^n)}.$$
Now, since $u,v\in B_R$, we get
$$ \n \Phi(u)-\Phi(v)\n_X \leq C_2 T R^{2k} \n u-v\n^{2}_{L^\infty_t([0,T])H^s_x(\mathbb{R}^n)}.$$
Finally, by suitably choosing the time $T$ we get that the operator $\Phi$ is a contraction and the result follows.


\section{The case $b\not\equiv 0$: Smoothing estimates}

We now consider the more general case \eqref{P2} where $b\not\equiv 0$. For simplicity of notation let us rename 
\begin{equation}
\label{PG}
\LL_\alpha=i\partial_t+t^\alpha\Delta_x+b(t,x)\cdot \nabla_x,
\end{equation}
where 
$\nabla_x=(\partial_{x_1},...,\partial_{x_n})$, $\Delta_x=\sum_{j=1}^n D_j^2$,  $b(t,x)=(b_1(t,x),...,b_n(t,x))$, $b_j\in C([0,T],C_b^\infty(\mathbb{R}^n))$ for all $j=1,...,n$, and $C^\infty_b(\mathbb{R}^n)=\{f\in C^\infty(\mathbb{R}^n): \partial^\alpha f\in L^\infty(\mathbb{R}^n),\,\forall\alpha\in\mathbb{Z}^n_+\}$. Moreover we assume that there exist $c_0>0$  and $\sigma>1$ such that, for all $j=1,...,n$,
$$|\mathsf{Im}\,\partial_x^\gamma b_j(t,x)|, |\mathsf{Re}\,\partial_x^\gamma b_j(t,x)|\lesssim t^\alpha \langle x \rangle^{-\sigma-|\gamma|}, \quad\text{for}\,\, |\gamma|=0,1.$$
 Our goal here is to prove some weighted smoothing estimates similar to the previous ones for the operator above. 
 
 \begin{remark}
 The smoothing estimates we are going to prove in this section are better than the ones proved in the case $b\equiv 0$, in the sense that for the non homogenous term we will be able to obtain some smoothing estimates with a gain of one derivative with respect to the regularity of initial data in the IVP. The smoothing estimates below can, of curse, be applied to the case $b\equiv 0$ as well. Moreover, one can prove these estimates in a more direct way in the case $b\equiv 0$. However, since we are interested in the case $b\not\equiv 0$ in which a direct proof is not applicable, we shall give the proof of the result for the general case directly. 
 \end{remark}

As we shall see, a key point in the proof of the smoothing properties for \eqref{PG} is the use of Doi's lemma (Lemma \ref{Doi} in the Appendix). 

We use Lemma \ref{Doi} on the symbol $a^w:=a=a_2+ia_1+a_0$ such that $a_2(x,\xi)=|\xi|^2$ and $a_1=a_0=0$.
In this case conditions (B1) and (B2) of Lemma \ref{Doi} are trivially satisfied, while (A6) holds with $q(x,\xi)=x\cdot\xi \langle \xi\rangle^{-1}$. Therefore, by Lemma \ref{Doi} with $\lambda'(|x|)= C' \langle x \rangle^{-\sigma}$ (see Remark \ref{remDoi}), with $C'$ is to be chosen later, we get that there exists $p\in S^0$ and $C>0$ such that \eqref{Doi1.1} holds.

We then consider the pseudo-differential operator $K$ with symbol $K(x,\xi)=e^{p(x,\xi)}\Lambda^s$, where $\Lambda^s:=\langle \xi\rangle^s$ and $p(x,\xi)$ is given by Doi's lemma, and define the norm $N$ on $H^s(\mathbb{R}^n)$, equivalent to the standard one (see \cite{KPRV} for the proof of the equivalence), given by
\begin{equation}
\label{normN}
N(u)^2=\n Ku\n^2_{0}+\n u\n_{s-1}^2,
\end{equation}
where $\n \cdot\n_{s}$ stands for the standard norm in the Sobolev space $H^{s}(\mathbb{R}^n)$.

Finally, following the technique used in \cite{KPRV}, we make use of the norm $N$ to prove smoothing properties of the solutions of the IVP \eqref{IVP2}
that we stated in Theorem \ref{thb}.

The proof of this theorem is essentially reduced to the proof of the following lemma.

\begin{lemma} \label{lemma}
Let $s\in\mathbb{R}$, $\lambda(|x|):=\langle x \rangle^{-\sigma}$ and $\sigma>1$ such that \eqref{condb} holds. Then there exists $C_1, C_2>0$ such that, for all $u\in C([0,T];H^{s+2}(\mathbb{R}^n))\bigcap C^1([0,T];H^s(\mathbb{R}^n))$, we have

  \begin{equation}
\label{sm1}
\sup_{0\leq t\leq T} \n u(t)\n_s\leq C_1 e^{C_2(\frac{T^{\alpha+1}}{\alpha+1}+T)}\left( \n u_0\n_s+\int_0^T \n (\partial_t-it^\alpha\Delta_x-ib(t,x)\cdot \nabla_x)u(t,\cdot)\n_s dt\right);
\end{equation}
  \begin{equation}
\label{sm1.1}
\sup_{0\leq t\leq T} \n u(t)\n_s\leq C_1 e^{C_2(\frac{T^{\alpha+1}}{\alpha+1}+T)}\left( \n u(\cdot,T)\n_s+\int_0^T \n (\partial_t-it^\alpha\Delta_x-ib(t,x)\cdot \nabla_x)^*u(t,\cdot)\n_s dt\right);
\end{equation}
 \begin{equation}
\label{sm2}
\sup_{0\leq t\leq T} \n u(t)\n^2_s +\int_0^T \int_{\mathbb{R}^n}  t^\alpha \left| \Lambda^{s+1/2} u \right|^2 \lambda(|x|)  dx\,dt
\end{equation}
$$\leq C_1 e^{C_2(\frac{T^{\alpha+1}}{\alpha+1}+T)}\left( \n u_0\n^2_s+\int_0^T \n (\partial_t-it^\alpha\Delta_x-ib(t,x)\cdot\nabla_x)u(t,\cdot)\n^2_s dt\right);$$
   \begin{equation}
\label{sm3}
 \sup_{0\leq t\leq T} \n u(t)\n^2_s +\int_0^T \int_{\mathbb{R}^n}  t^\alpha \left| \Lambda^{s+1/2} u \right|^2 \lambda(|x|)dx\,dt
 \end{equation}
$$\leq C_1 e^{C_2\frac{T^{\alpha+1}}{\alpha+1}}\left( \n u_0\n^2_s+\int_0^T \int_{\mathbb{R}^n} t^{-\alpha} \lambda(|x|)^{-1}   \left| \Lambda^{s-1/2} (\partial_t-it^\alpha\Delta_x-ib(t,x)\cdot \nabla_xu(t, \cdot)\right|^2 dx\, dt\right).$$

\end{lemma}

\proof
Recall that we defined the norm $N(u)$ on $H^s(\mathbb{R}^n)$ to be the one given in \eqref{normN}.
We then estimate the quantity 
$$\partial_t N(u)^2= \partial_t\n Ku\n^2_{0}+\partial_t\n u\n_{s-1}^2=I+II.$$
 For the term $II$ we have
 $$II=\partial_t \n u\n_{s-1}^2=\partial_t\la \Lambda^{s-1}u,\Lambda^{s-1}u\ra=2\mathsf{Re}\la \Lambda^{s-1}\partial_t u, \Lambda^{s-1}u\ra$$
$$=2\mathsf{Re}\la \Lambda^{s-1}(it^\alpha\Delta_x+ib(t,x)\cdot \nabla_x u), \Lambda^{s-1}u\ra$$
$$= 2\underbrace{\mathsf{Re}\la it^\alpha \Lambda^{s-1} \Delta_x u, \Lambda^{s-1}u\ra}_{=0}
-2\mathsf{Re}\la \Lambda^{s-1}b(t,x)\cdot D_x u, \Lambda^{s-1} u\ra+2\mathsf{Re}\la \Lambda^{s-1} f,\Lambda^{s-1} u\ra$$
$$\leq C t^\alpha \n u\n^2_s+2\mathsf{Re}\la \Lambda^{s-1} f,\Lambda^{s-1} u\ra,$$
where (recall) $D_x=(D_{x_1},...,D_{x_n})=(-i\partial_{x_1},...,-i\partial_{x_n})$ and 
$b(t,x)\cdot D_x=\sum_{j=1}^nb_j(t,x)D_{x_j}$, with $b_j\in C([0,T];C^\infty_b(\mathbb{R}^n))$ such that \eqref{condb} holds.

Observe that the following estimates hold:
\begin{equation}
\label{estf1}
 2\mathsf{Re}\la \Lambda^{s-1} f,\Lambda^{s-1} u\ra\leq 2 \n f\n_{s-1}\n u\n_{s-1}\leq C N(f)N(u),
\end{equation}
and
$$ 2\mathsf{Re}\la \Lambda^{s-1} f,\Lambda^{s-1} u\ra= 2\mathsf{Re}\la t^{-\alpha/2}\lambda(|x|)^{-1/2}\Lambda^{s-1/2} f,t^{\alpha/2}\lambda(|x|)^{1/2}\Lambda^{s-3/2} u\ra$$
$$\leq \n t^{-\alpha/2}\lambda(|x|)^{-1/2}\Lambda^{s-1/2} f\n_0^2+\n t^{\alpha/2}\lambda(|x|)^{1/2}\Lambda^{s-3/2} u\n_0^2$$
\begin{equation}
\label{estf2}
\leq \la t^{-\alpha}\lambda(|x|)^{-1}\Lambda^{s-1/2} f,\Lambda^{s-1/2}f\ra +t^\alpha N( u)^2.
\end{equation}
Therefore, by using \eqref{estf1} and \eqref{estf2}, we get that $II$ can be estimated  by 

\begin{equation}
\label{estII2}
II\leq C t^\alpha N(u)^2+C'\min \{N(f)N(u);  \la t^{-\alpha}\lambda(|x|)^{-1}\Lambda^{s-1/2} f,\Lambda^{s-1/2}f\ra \},
\end{equation}
with $C$ and $C'$ new suitable constants.

We now consider the term $I$

\begin{equation}
\label{dtK}
\partial_t\n Ku\n ^2_0=2\mathsf{Re}\la \partial_t Ku,Ku\ra=2\mathsf{Re}\la K\partial_t u,Ku\ra
\end{equation}
$$=2\mathsf{Re}\la K(it^\alpha \Delta_x+ib(t,x)\cdot \nabla_x)u) ,Ku\ra+2\mathsf{Re}\la Kf,Ku\ra$$
$$=2\mathsf{Re}\la it^\alpha[K,\Delta_x]u,Ku\ra+\underset{=0}{\underbrace{2\mathsf{Re}\la it^\alpha \Delta_x Ku,Ku\ra}}$$
$$-2\mathsf{Re}\la K\,b(t,x)\cdot D_xu,Ku\ra+2\mathsf{Re}\la Kf,Ku\ra$$
$$=2\mathsf{Re}\la it^\alpha[K,\Delta_x]u,Ku\ra-2\mathsf{Re}\la [K,b(t,x)\cdot D_x]u,Ku\ra$$
$$-2\mathsf{Re}\la b(t,x)\cdot D_x\,Ku,Ku\ra+2\mathsf{Re}\la Kf,Ku\ra,$$
 \\
and estimate the term $2\mathsf{Re}\la [K,b(t,x)\cdot D_x]u,Ku\ra$ in the the fifth line of \eqref{dtK}.
Recall that, given two symbols $p_1\in S^{m_1}, p_2\in S^{m_2}$ associated with two operators $P_1$ and $P_2$, then we have that the symbol of the commutator $[P_1,P_2](x,D)$ is given by $-i\{p_1,p_2\}(x,\xi)+p_3(x,\xi)$, where $p_3\in S^{m_1+m_2-2}$. Therefore, since $K(x,\xi)=e^{p(x,\xi)}(1+|\xi|^2)^{s/2}$ and the symbol of $b\cdot D_x=b(t,x)\cdot D_x$ is $\sum_{j=1}^nb_j(t,x)\xi_j$, we have that the operator $[K,b(t,x)\cdot D_x]$ is of order $s$ and has symbol 

\begin{equation}\label{reminder}
[K,b(t,x)\cdot D_x](t,x,\xi)=-i\{K(x,\xi),\sum_{j=1}^nb_j(t,x)\xi_j \}+ r_{s-1}(t,x,\xi)
\end{equation}
$$=-i\sum_{k=1}^n\Big[ e^{p(x,\xi)} (\Lambda^s(\xi)\partial_{\xi_k}p(x,\xi)+s \Lambda^{s-1}(\xi)\xi_k) \sum_{j=1}^n(\partial_{x_k}b_j(t,x))\xi_j$$
$$- e^{p(x,\xi)}(\partial_{x_k}p(x,\xi))\Lambda^s(\xi)b_k(t,x)\Big]+r_{s-1}(t,x,\xi).$$

Therefore, by the properties of $b(t,x)$ (recall that $b\in C^\infty_b$ and is bounded, together with its derivatives in space, by $t^\alpha \lambda(|x|)$ we get
$$-2\mathsf{Re}\la [K,b(t,x)D_x]u,Ku\ra\leq C t^\alpha \n u\n_s^2,$$
 where we used $\n r_{s-1}(t,x,D)u\n_0\leq Ct^\alpha \n u\n_{s-1}$ (this estimate is deduced by using the properties of $b$ following \cite{L} Theorem 1.1.20 pag.14).
 
Note also that, once more by using the pseudo-differential calculus, we get $[K,\Delta_x](x,D)=[p,\Delta_x]K(x,D)+r_{s}(x,D)$, where $r_s$ is of order $s$ and $p=p(x,D)$ is the operator of order $0$ appearing in the definition of the norm $N(\cdot)$.
 
Now we can estimate \eqref{dtK} in the following way

$$\eqref{dtK}\leq C t^\alpha \n u\n_s^2+2\mathsf{Re}\la(it^\alpha [p,\Delta_x](x,D)-b(t,x)\cdot D_x)Ku,Ku\ra + |2\mathsf{Re}\la it^\alpha r_s(x,D)u,Ku\ra|$$
\begin{equation}
\label{dtK2}
 \leq C t^\alpha \n u\n_s^2+2\mathsf{Re}\la(it^\alpha [p,\Delta_x](x,D)-b(t,x)D_x)Ku,Ku\ra, 
\end{equation}
where $C$ is a new suitable positive constant.

We denote $Q(x;D):=it^\alpha [p,\Delta_x](x,D)-b(t,x)\cdot D_x$ whose symbol is such that

$$\mathsf{Re}\,Q(x,\xi)=\mathsf{Re}\left(i t^\alpha (-i)\{ p, -|\xi|^2\}(x,\xi)-b(t,x)\cdot \xi\right) +r_0$$
$$\leq -t^\alpha \{ p, |\xi|^2\}(x,\xi)+|\mathsf{Re}\,b(t,x)\cdot \xi | +r_0$$
$$\leq -t^\alpha \{ p, |\xi|^2\}(x,\xi)+|\mathsf{Re}\, b(t,x)||\xi | +C_4$$
$$\underset{\text{by \eqref{Doi1.1} }}{\leq } -C't^\alpha \lambda(|x|)|\xi| +C_2t^\alpha +C_0t^\alpha\lambda(|x|)|\xi|+C$$
$$\leq t^\alpha (C_0-C') \lambda(|x|)|\xi|+C_2t^\alpha  +C_4$$
$$\leq -C t^\alpha \lambda(|x|)|\xi| +C_2t^\alpha +C $$
$$\leq -Ct^\alpha \lambda(|x|) (1+|\xi|) +C t^\alpha \lambda(|x|)+C_2t^\alpha +C_4 $$
$$\leq -Ct^\alpha \lambda(|x|) (1+|\xi|^2)^{1/2} +C_3t^\alpha +C_4 $$
$$=t^\alpha(-C\lambda(|x|)(1+|\xi|^2)^{1/2}+C_3)+C_4$$
where we chose $C'$ (which is possible by Doi's lemma, see Remark \ref{remDoi}) in order to have $C_0-C' <0$.

Due to the property of the symbol of $Q$ we can apply the G\aa rding inequality and get

$$2\mathsf{Re} \la Q(x,D)Ku,Ku \ra\leq -Ct^\alpha \la \lambda(|x|)\Lambda^1 Ku, Ku\ra +C_3t^\alpha \n Ku\n^2_0+C_4\n Ku\n_0^2$$
$$\leq -Ct^\alpha \la \lambda(|x|)\Lambda^1 Ku, Ku\ra +C_3t^\alpha\n u\n_s^2+C_4\n u\n_s^2$$

Since $\lambda\in C^\infty_b$, by using the symbolic calculus we get that 
$$ \lambda(|x|)\Lambda^1(x,D)= (\lambda(|x|)^{1/2}\Lambda^{1/2})^{2}(x,D)+ r_0(x,D),$$

where $r_0(x,D)$ has order 0.

Then, by the latter property, we get
\begin{equation}
\label{estQ}
2\mathsf{Re} \la Q(x,D)Ku,Ku \ra\leq 
-Ct^\alpha \n \lambda(|x|)^{1/2}\Lambda^{1/2} Ku\n_0^2+ C_3t^\alpha\n u\n_s^2
+C_4\n u\n_s^2,
\end{equation}
where $C>0$ is a new suitable constant.

By plugging \eqref{estQ} in \eqref{dtK2} we get
\begin{equation}
\label{ }
\partial_t\n Ku\n_0\leq 
C t^\alpha N(u)^2+ C'N(u)^2-C''t^\alpha \n \lambda(|x|)^{1/2}\Lambda^{1/2} Ku\n_0^2+C'''N(f)N(u),
\end{equation}
where in the latter we used $\mathsf{Re}\langle Kf,Ku\ra\leq C'''N(f)N(u)$.

Finally, by using \eqref{estII2} and the equivalence between the norms $\n \cdot\n_s$ and $N(\cdot)$ (see \cite{KPRV} pag.390), we obtain
\begin{equation}
\label{estN}
\partial_t N(u)^2=\partial_t \n Ku\n^2+\partial_t \n u\n_{s-1}^2
\end{equation}
$$\leq C t^\alpha N(u)^2+ C'N(u)^2-C''t^\alpha \n \lambda(|x|)^{1/2}\Lambda^{1/2} Ku\n_0^2+C'''N(f)N(u)+ $$
$$+C_3\min \{N(f)N(u);  \la t^{-\alpha}\lambda(|x|)^{-1}\Lambda^{s-1/2} f,\Lambda^{s-1/2}f\ra \},$$
where the constants are (eventually) new suitable constants.

From \eqref{estN} we will get \eqref{sm1}, \eqref{sm1.1} and \eqref{sm2} as we shall prove below.

\proof[Proof of \eqref{sm1}]
As regards \eqref{sm1} we observe that, from \eqref{estN} we have
$$\partial_t N(u)^2\leq C_1(t^\alpha +1) N(u)^2+  C_2N(u)N(f)$$
(again $C_1$ and $C_2$ new constants) which gives, 
$$2\partial_t N(u)\leq C_1(t^\alpha +1) N(u)+ C_2N(f)$$
and
$$\partial_t\left( 2e^{-\frac1 2 C_1(t^{\alpha+1}/(\alpha+1)+t)}N(u)\right)\leq C_2e^{-\frac1 2 C_1(t^{\alpha+1}/(\alpha+1)+t)}  N(f).$$
Hence, by integrating in time from 0 to $t$ we get

$$N(u(t))\leq C e^{\frac1 2 C_1(t^{\alpha+1}/(\alpha+1)+t)}\left[ N(u(0))+C_2\int_0^t e^{-\frac1 2 C_1(s^{\alpha+1}/(\alpha+1)+s)}N(f)ds \right]$$
$$\leq C' e^{\frac1 2 C_1(t^{\alpha+1}/(\alpha+1)+t)}\left[ N(u(0))+\int_0^t N(f)ds \right],$$

which finally gives \eqref{sm1} by the equivalence of the norms.
\endproof
\proof[Proof of \eqref{sm1.1} ]
The proof of \eqref{sm1.1}  follows from \eqref{sm1}  applied to the adjoint operator and with $u(t,\cdot)$ replaced by $u(T-t,\cdot)$.
\endproof

\proof[proof of \eqref{sm2}]
To obtain \eqref{sm2} we first observe that there exists a pseudo-differential operator $\tilde{K}$ such that
$$ I=\tilde{K}K+\Psi_{r_{-1}},$$
where $\Psi_{r_{-1}}$ is an operator with symbol $r_{-1}$ of order $-1$ (see \cite{KPRV} pag.390 for the proof).
By using this property we get

\begin{equation}
\label{Ks}
\n \lambda(|x|)^{1/2}\Lambda^{s+1/2} u\n_0\leq \n ( \lambda(|x|)^{1/2}\Lambda^{1/2})(\Lambda^s \tilde{K})(K\Lambda^{1/2})u\n_0
+ O(N(u))
\end{equation}
$$\leq \n(\Lambda^s \tilde{K})( \lambda(|x|)^{1/2}\Lambda^{1/2})(K\Lambda^{1/2})u\n_0+c N(u)$$
$$\leq c  \left(\n( \lambda(|x|)^{1/2}\Lambda^{1/2})(K\Lambda^{1/2})u\n_0+ N(u)\right),$$
where, in the second line, we used the fact that $[\Lambda^s \tilde{K}, \lambda(|x|)^{1/2}\Lambda^{1/2}]K\Lambda^{1/2}$ is a pseudo-differential operator of order $s$ together with the equivalence of the norms $\n\cdot\n_s$ and $N(\cdot)$.
Therefore, again from \eqref{estN}, we have
$$\partial_t N(u)^2\leq C_1(t^\alpha +1) N(u)^2-C''t^\alpha \n \lambda(|x|)^{1/2}\Lambda^{1/2} Ku\n_0^2+C'''N(f)N(u)$$
$$\leq C_1(t^\alpha +1) N(u)^2-C_2t^\alpha \la \lambda(|x|)^{1/2}\Lambda^{s+1/2} u, \lambda(|x|)^{1/2}\Lambda^{s+1/2} u\ra +C_3N(u)^2+C_4N(f)^2$$
$$\leq C_1(t^\alpha +1) N(u)^2-C_2t^\alpha \la \lambda(|x|)^{1/2}\Lambda^{s+1/2} u, \lambda(|x|)^{1/2}\Lambda^{s+1/2} u\ra +C_4N(f)^2,$$
where the constants are new suitable constant. Hence

$$\partial_t N(u)^2+C_2 \la t^{\alpha/2}\lambda(|x|)^{1/2}\Lambda^{s+1/2} u,t^{\alpha/2} \lambda(|x|)^{1/2}\Lambda^{s+1/2} u\ra \leq C_1(t^\alpha +1) N(u)^2+C_4N(f)^2,$$
so that by integrating in time from 0 to $t$,

$$N(u(t))^2+ C_2 e^{\frac1 2 C_1(t^{\alpha+1}/(\alpha+1)+t)} \times$$
\begin{equation}
\label{iii}
\times \int_0^t e^{-\frac1 2 C_1(s^{\alpha+1}/(\alpha+1)+s)} \la s^{\alpha/2}\lambda(|x|)^{1/2}\Lambda^{s+1/2} u,s^{\alpha/2} \lambda(|x|)^{1/2}\Lambda^{s+1/2} u\ra ds
\end{equation}
$$\lesssim e^{\frac1 2 C_1(t^{\alpha+1}/(\alpha+1)+t)}\left[ N(u(0))^2+\int_0^t e^{-\frac1 2 C_1(s^{\alpha+1}/(\alpha+1)+s)}N(f)^2ds \right]$$
$$\lesssim  e^{\frac1 2 C_1(T^{\alpha+1}/(\alpha+1)+T)}\left[ N(u(0))^2+\int_0^t N(f)^2ds \right].$$
From the previous estimate we get
\begin{equation}
\label{iii1}
\sup_{t\in[0,T]}N(u(t))^2 \lesssim e^{\frac1 2 C_1(T^{\alpha+1}/(\alpha+1)+T)}\left[ N(u(0))^2+\int_0^T N(f)^2ds \right].
\end{equation}
Moreover, the second term on the LHS of \eqref{iii} satisfies
$$e^{\frac1 2 C_1(t^{\alpha+1}/(\alpha+1)+t)} \int_0^t e^{-\frac1 2 C_1(s^{\alpha+1}/(\alpha+1)+s)} \la s^{\alpha/2}\lambda(|x|)^{1/2}\Lambda^{s+1/2} u,s^{\alpha/2} \lambda(|x|)^{1/2}\Lambda^{s+1/2} u\ra ds$$
$$\geq e^{\frac1 2 C_1(t^{\alpha+1}/(\alpha+1)+t)}(\inf_{s\in[0,T]]} e^{-\frac1 2 C_1(s^{\alpha+1}/(\alpha+1)+s)} )\times$$
$$\times \int_0^t \la s^{\alpha/2}\lambda(|x|)^{1/2}\Lambda^{s+1/2} u,s^{\alpha/2} \lambda(|x|)^{1/2}\Lambda^{s+1/2} u\ra ds$$
$$\geq \int_0^t \la s^{\alpha/2}\lambda(|x|)^{1/2}\Lambda^{s+1/2} u,s^{\alpha/2} \lambda(|x|)^{1/2}\Lambda^{s+1/2} u\ra ds.$$
Therefore, using the previous inequality and \eqref{iii},
\begin{equation}
\label{iii2}
C_2\int_{0}^T \la s^{\alpha/2}\lambda(|x|)^{1/2}\Lambda^{s+1/2} u,s^{\alpha/2} \lambda(|x|)^{1/2}\Lambda^{s+1/2} u\ra ds
\end{equation}
$$=C_2 \sup_{t\in[0,T]} \int_0^t \la s^{\alpha/2}\lambda(|x|)^{1/2}\Lambda^{s+1/2} u,s^{\alpha/2} \lambda(|x|)^{1/2}\Lambda^{s+1/2} u\ra ds$$
$$\leq C_2 \sup_{t\in[0,T]}  e^{\frac1 2 C_1(t^{\alpha+1}/(\alpha+1)+t)} \int_0^t e^{-\frac1 2 C_1(s^{\alpha+1}/(\alpha+1)+s)} \la s^{\alpha/2}\lambda(|x|)^{1/2}\Lambda^{s+1/2} u,s^{\alpha/2} \lambda(|x|)^{1/2}\Lambda^{s+1/2} u\ra ds$$
$$\underset{\eqref{iii}}{\leq} e^{\frac1 2 C_1(T^{\alpha+1}/(\alpha+1)+T)}\left[ N(u(0))+\int_0^T N(f)^2ds \right].$$
Finally, by summing up \eqref{iii1} and \eqref{iii2} we get \eqref{sm2}.
\endproof

\proof[Proof of \eqref{sm3}]
To prove \eqref{sm3}, denoting by $\lambda:=\lambda(|x|)$, we write
\begin{equation}
\label{iv1}
2\mathsf{Re}\la Kf, Ku\ra=2\mathsf{Re}\la t^{\alpha/2}\lambda^{1/2}\Lambda^{1/2}Kf,t^{-\alpha/2}\lambda^{-1/2}\Lambda^{-1/2} Ku\ra
\end{equation}
$$\leq \varepsilon \n  t^{\alpha/2}\lambda^{1/2}\Lambda^{1/2}Ku \n_0^2+\frac 1 \varepsilon \n t^{-\alpha/2}\lambda^{-1/2}\Lambda^{-1/2} Kf\n_0^2$$
$$= \varepsilon \n  t^{\alpha/2}\lambda^{1/2}\Lambda^{1/2}K\Lambda^{-s-1/2}\Lambda^{s+1/2}u\n_0^2$$
$$+\frac 1 \varepsilon \n t^{-\alpha/2}\lambda^{-1/2}\Lambda^{-1/2} K\Lambda^{-s+1/2}\Lambda^{s-1/2}f\n_0^2,$$

Since $\Lambda^{1/2}K\Lambda^{-s-1/2}$ and $\Lambda^{-1/2} K\Lambda^{-s+1/2}$ are both pseudo-differential operators of order 0 in $x$, we have 
$$ t^{\alpha/2}\lambda^{1/2}\Lambda^{1/2}K\Lambda^{-s-1/2}= \Lambda^{1/2}K\Lambda^{-s-1/2}t^{\alpha/2}\lambda^{1/2} +t^{\alpha/2}\Psi_{r_{-1}},$$
where $\Psi_{r_{-1}}$ denotes an operator of order $-1$ in the space variable. Of course the same property holds for the operator $t^{-\alpha/2}\lambda^{-1/2}\Lambda^{-1/2} K\Lambda^{-s+1/2}$.

We use these properties in \eqref{iv1} to get 

\begin{equation}
\label{iv2}
2\mathsf{Re}\la Kf, Ku\ra=2\mathsf{Re}\la t^{\alpha/2}\lambda^{1/2}\Lambda^{1/2}Kf,t^{-\alpha/2}\lambda^{-1/2}\Lambda^{-1/2} Ku\ra
\end{equation}
$$\leq c_1 \varepsilon \n t^{\alpha/2}\lambda^{1/2}\Lambda^{s+1/2}u\n_0^2+ c_2\frac 1 \varepsilon \n t^{-\alpha/2}\lambda^{-1/2} \Lambda^{s-1/2}f\n_0^2 +c_3 t^\alpha \n u\n^2_s.$$

By using \eqref{Ks} and \eqref{iv2} in \eqref{estN}, and the equivalence between the norms $N(\cdot)$ and $\n \cdot \n_s$, we obtain
$$\partial_tN(u)^2+ (c_0-c_1\varepsilon) \n t^{\alpha/2}\lambda^{1/2}\Lambda^{s+1/2}u\n_0^2 
\leq c_3 t^\alpha N( u)^2 +c_2\frac 1 \varepsilon \n t^{-\alpha/2}\lambda^{-1/2} \Lambda^{s-1/2}f\n_0^2,$$
where $c_j$, $j=0,1,2,3$ are new suitable constants, and we choose $\varepsilon >0$ such that $c_0-c_1\varepsilon\geq c>0$.

Since $\varepsilon$ is now fixed, we have
$$\partial_tN(u)^2+ c \n t^{\alpha/2}\lambda^{1/2}\Lambda^{s+1/2}u\n_0^2 \leq c_3 t^\alpha N( u)^2 +c_2\frac 1 \varepsilon \n t^{-\alpha/2}\lambda^{-1/2} \Lambda^{s-1/2}f\n_0^2,$$
which gives, by integrating in time from $0$ to $t$, and by using the same argument as in the proof of \eqref{sm2}, the proof of \eqref{sm3}.
\endproof
The proof is then complete.
\endproof

Lemma \ref{lemma} allows to prove the  well-posedness and smoothing result in Theorem \ref{thb}.

\proof[Proof of Theorem \ref{thb}]
From \eqref{sm1} of Lemma \ref{lemma} we immediately get the uniqueness of the solution. In fact, let $u$ be a solution of the homogeneous IVP for \eqref{PG}, i.e., with $f=0$, and initial data $u_0=0$. Then, by  \eqref{sm1} of Lemma \ref{lemma}, we get $u=0$ and thus the uniqueness (even in the general case $f\neq0$ and $u_0\neq 0$).

About the existence, we get the results by using density arguments as we will prove below.\\

\vspace{0.4cm}
\textit{Case 1:} $f\in \mathcal{S}(\mathbb{R}^{n+1})$ and $u_0\in \mathcal{S}(\mathbb{R}^{n})$.\\	
We consider the subspace
$$E=\{P^*\varphi;\,\varphi\in C_0^\infty(\mathbb{R}^{n}\times[0,T))\}=(\partial_t-it^\alpha\Delta_x+b(t,x)\cdot D_x)^*(C_0^\infty(\mathbb{R}^{n+1}))\subset L^1([0,T];H^{-s}(\mathbb{R}^n)$$
and the linear functional
$$\ell^*: E\rightarrow \mathbb{C},\quad \ell^*(P^*\varphi)=\int_0^T \la f,\varphi\ra_{L^2\times L^2} dt+\la u_0,\varphi(\cdot,0)\ra_{L^2\times L^2}.$$

Then, by \eqref{sm1.1} of Lemma \ref{lemma} (applied on $\varphi$) with $s$ replaced by $-s$, for $\eta=P^*\varphi$, with $\varphi \in C_0^\infty(\mathbb{R}^{n}\times[0,T))$ we get
$$|\ell^*(\eta)|\leq \n f\n_{(L^1[0,T];H^s_x)} \sup_{t\in[0,T]}\n\varphi\n_{H_x^{-s}}+\n u_0\n_{H_x^s}\n\varphi(0)\n_{H_x^{-s}}$$
$$\leq e^{C(T^{\alpha+1}/(\alpha+1)+T)}\left(\n f\n_{L^1_t([0,T];H^s_x)}+\n u_0\n_{H_x^s}\right)\n\eta\n_{L^1_t([0,T];H^{-s}_x)},$$
which gives the continuity of $\ell^*$ on $E$ (the last inequality follows from  \eqref{sm1.1} of Lemma \ref{lemma} applied both on the term $\sup_{t\in[0,T]}\n\varphi\n_{H_x^{-s}}$ and $\n\varphi(0)\n_{H_x^{-s}}$ together with the compactness of the support of $\varphi$). By the Hahn-Banach theorem we can extend $\ell^*$ on $L^1([0,T]:H^{-s}(\mathbb{R}^n))$ and finally get the existence of $u\in L^1([0,T];H^{-s}(\mathbb{R}^n))^*=L^\infty([0,T];H^{s}(\mathbb{R}^n))$ such that
$$\ell^*(P^*\varphi)=\la u, P^*\varphi\ra_{L^2\times L^2}=\int_0^T \la f,\varphi\ra_{L^2\times L^2} dt+\la u_0,\varphi(\cdot,0)\ra_{L^2\times L^2},\quad \forall\varphi\in C_0^\infty(\mathbb{R}^{n}\times[0,T))$$
and thus $Pu=f$ in the sense of distributions for $0<t<T$.

Notice that $Pu\overset{\mathcal{D}'}{=}f$  means that $(\partial_t-it^\alpha\Delta_x+b(t,x)\cdot D_x)u\overset{\mathcal{D}'}{=}f$ (as distributions on $C_0^\infty([0,T]\times\mathbb{R}^n)$), therefore, since $f\in \mathcal{S}(\mathbb{R}^{n+1})$, we have that $\partial_t u\in (L^\infty [0,T):H^{s-2}(\mathbb{R}^{n}))$, which gives $u \in (C([0,T):H^{s-2}(\mathbb{R}^{n}))$. We then use the equation once more, that is $\partial_tu=it^\alpha\Delta_x+b(t,x)\cdot D_xu +f$, and get, by doing the same consideration, that $u\in (C^1[0,T):H^{s-4}(\mathbb{R}^{n}))$ and $u(x,0)=u_0(x)$. Finally, since $u_0\in H^s(\mathbb{R}^{n})$, repeating the previous argument with $s+4$ in place of $s$ we conclude that there exists a solution $u$ of the IVP associated to \eqref{PG} to which parts (i)-(iv) of Lemma \ref{lemma} apply.

\vspace{0.4cm}
\textit{Case 2:} $f\in L^1([0,T];H^s(\mathbb{R}^{n}))$ and $u_0\in H^s(\mathbb{R}^{n})$.\\	
In this case we take two sequences $f_j\in \mathcal{S}(\mathbb{R}^{n+1})$, $v_j\in \mathcal{S}(\mathbb{R}^{n})$ such that $f_j\rightarrow f$ in $(L^1([0,T]):H^s(\mathbb{R}^{n})$ and 
$v_j\rightarrow u_0$ in $H^s(\mathbb{R}^{n})$.

By the arguments of case 1 we find a solution $u_j$ of the IVP associated with \eqref{PG} with $f_j$ and $v_j$ in place of $f$ and $u_0$ respectively. Since $u_j$ satisfies  \eqref{sm1} of Lemma \ref{lemma}, we have that $u_j$ is a Cauchy sequence, therefore, passing to the limit, we get that $u=\lim_{j\rightarrow}u_j$ is a solution of the IVP with initial data $f$ and in initial data $u_0$ satisfying \eqref{sm1.1} of Lemma \ref{lemma}, which proves (ii) of the theorem.

\vspace{0.4cm}
\textit{Case 3:} $f\in L^2([0,T];H^s(\mathbb{R}^{n}))$ and $u_0\in H^s(\mathbb{R}^{n})$.\\	
Here we proceed as in case 2 where, instead, $f_j\in \mathcal{S}(\mathbb{R}^{n+1})$ is such that $f_j\rightarrow f$ in $(L^2([0,T]);H^s(\mathbb{R}^{n})$. Under this hypothesis we get point (ii) of the theorem, that is, it exists a solution $u\in (C[0,T):H^{s}(\mathbb{R}^{n}))$ satisfying \eqref{sm2} of Lemma \ref{lemma}.

\vspace{0.4cm}
\textit{Case 4:} $\Lambda^{s-1/2}f\in (L^2(\mathbb{R}^n\times [0,T]): t^{-\alpha}\lambda(|x|)^{-1}dxdt)$ and $u_0\in H^s(\mathbb{R}^{n})$.\\	
In this case it is possible to prove that there exists $g_j\in\mathcal{S}(\mathbb{R}^{n+1})$ such that $g_j\rightarrow \Lambda^{s-1/2}f$ in $(L^2(\mathbb{R}^n)\times [0,T]: t^{-\alpha}\lambda(|x|)^{-1}dxdt)$. Applying once again the procedure used in case 1 with $f_j$ replaced by $\Lambda^{-s+1/2} g_j$ in  \eqref{sm3} of Lemma \ref{lemma}, and passing to the limit, we finally get point (iii) of Theorem \ref{thb}.

\endproof

\section{The case $b\neq 0$: local well-posedness  of the nonlinear Cauchy problem}
We now analyze the local well-posedness of the IVP 
\begin{equation}
\label{IVP3}
\left\{ \begin{array}{ll}
\partial_tu=it^\alpha\Delta_xu +ib(t,x)\cdot \nabla_xu+P(u,\bar{u},\nabla u, \nabla \bar{u})& \\
u(0,x)=u_0(x). & 
\end{array}\right .
\end{equation}
under the previous hypotheses on the term $b$ in \eqref{IVP3}. Note that, with some abuse of notation, the quantity $\nabla u$ in the nonlinear term is $\nabla u:=\mathsf{div}(u):=\sum_{j=1}^n\partial_{x_j}u$, and, simlarily, the quantity $\nabla \bar{u}$.

First we shall consider the case $P(u,\bar{u},\nabla u, \nabla \bar{u})=P(u)=\pm u|u|^{2k}$, $k\geq1$, which is treated in Theorem \ref{teo1} and, afterwards, the case $P(u,\bar{u},\nabla u, \nabla \bar{u})=\pm t^\beta \nabla_xu \cdot u^{2k}$ in Theorem \ref{teo2}.
The proof is based on the contraction argument, which, once again, is obtained through the use of the smoothing estimates proved in the previous section. 

\begin{remark}
Observe that in the proof below we will assume that the solution of the homogeneous problem is again of the form $W_\alpha(t)u_0$. Since we know from Theorem \ref{thb} that the solution of the linear problem exists, we assume that there exists a two-parameter family of operators, denoted $W_\alpha(t,\tau)$, giving the solution at time $t$ of the homogeneous problem with initial condition at time $\tau$ (recall, $W_\alpha(t,0):=W_\alpha(t)$). Under this assumption one can prove that Duhamel's formula still applies, therefore it makes sense to consider the operator $\Phi$ in the form given above. This is important to keep in mind since we will apply the same strategy in the subsequent case as well, that is, in the case  $P(u,\bar{u},\nabla u, \nabla \bar{u})=\pm t^\beta \nabla_xu \cdot u^{2k}$.
\end{remark}

\proof[Proof of Theorem \ref{teo1}]
We shall make use of the result in Theorem \ref{thb} concerning the linear problem to prove the result in the nonlinear case. Once again, we give the proof in the defocusing case since the proof in the focusing case applies with no modifications.
According to Theorem \ref{thb} we have the local well-posedness in $H^s$, $s>n/2$, for the linear IVP \eqref{IVP2}.
We now write the solution of \eqref{IVP3} as 
\begin{equation}
\label{Nu}
u(t,x)=W_\alpha(t)u_0+\int_0^tW_\alpha(t,\tau)P(u,\bar{u},\nabla u, \nabla \bar{u})d\tau,
\end{equation}
where $W_\alpha(t,\tau)$ is a new suitable two-parameter family of unitary operators.

Because of the previous assumption, solving the IVP \eqref{IVP3} with $P(u,\bar{u},\nabla u, \nabla \bar{u})=u|u|^{2k}$ is equivalent to find the solution of the integral equation

$$u(t,x)=W_\alpha(t)u_0+\int_0^tW_\alpha(t,\tau)u|u|^{2k}d\tau,$$
therefore, as in the proof of Theorem \ref{thb}, we look for the solution given by the fixed point of the map 
$$\Phi(u):= W_\alpha(t)u_0+\int_0^tW_\alpha(t,\tau)u|u|^{2k}d\tau,$$
defined on 
$$X_T^s:=\{u: [0,T]\times \mathbb{R}^n\rightarrow \mathbb{C}; \n u\n_{L^\infty_t H^s_x}<\infty,\,\left( \int_0^T\int_{\mathbb{R}^n}t^\alpha\lambda(|x|)|\Lambda^{s+1/2}u|^2 dx\,  dt \right)^{1/2}<\infty\},$$
where, recall, $\lambda(|x|):=\la x\ra ^\sigma$, with $\sigma>1$ and such that \eqref{condb} holds.

Proving the existence of a fixed point for $\Phi$ is, once more, equivalent to prove that the map $\Phi$ is a contraction, and, in order to do that, we will first show that $\Phi$ sends $X_T^s$ into itself.
The key point here will be to reduce ourself to a linear case to which the linear smoothing estimates apply.

Observe that, denoting by $v:=\Phi(u)$, we have that $v$ solves the linear problem

\begin{equation*}
\left\{ \begin{array}{ll}
\partial_tv-it^\alpha\Delta_xv -ib(t,x)\cdot \nabla_xv=u|u|^{2k}& \\
v(0,x)=u_0(x), & 
\end{array}\right .
\end{equation*}
and we can estimate $\n \Phi(u)\n_{X^s_T}=\n v\n_{X^s_T}$ through the smoothing estimates given by point (i) and (ii) of Theorem \ref{thb}.

In particular we define 
$$\n v\n_{X^s_T}:= \n v\n_{L^\infty_t H^s_x}+\left( \int_0^T\int_{\mathbb{R}^n}t^\alpha\lambda(|x|)|\Lambda^{s+1/2}v|^2 dx\,  dt\right)^{1/2} $$
$$=I+II,$$
where, by (i) of Theorem \ref{thb},
$$I= \n v\n_{L^\infty_t H^s_x}\leq Ce^{C' \frac{T^{\alpha+1}}{\alpha+1}+T}\left( \n u_0\n_s +\int_0^T \n u|u|^{2k}\n_s dt\right)$$
and, by (ii) of  Theorem \ref{thb},
$$II= \left(\int_0^T\int_{\mathbb{R}^n}t^\alpha\lambda(|x|)|\Lambda^{s+1/2}v|^2 dx\,  dt\right)^{1/2} $$
$$\leq Ce^{\frac{C'}{2} \left(\frac{T^{\alpha+1}}{\alpha+1}+T\right)}\left( \n u_0\n^2_s +\int_0^T \n u|u|^{2k}\n^2_s dt\right)^{1/2}.$$
$$\leq Ce^{C' \left(\frac{T^{\alpha+1}}{\alpha+1}+T\right)}\left( \n u_0\n^2_s +\left(\int_0^T \n u|u|^{2k}\n^2_s dt\right)^{1/2}\right).$$
Following the same computations of Theorem \ref{teo2} we have
$$\int_0^T \n u|u|^{2k}\n_{H^s_x} dt= \n u|u|^{2k}\n_{L^1_tH^s_x}\leq CT\n u\n_{L^\infty_t H^s_x}^{2k+1}$$
and
$$\int_0^T \n u|u|^{2k}\n^2_{H^s_x} dt\leq \n u|u|^{2k}\n_{L^\infty_tH^s_x}\n u|u|^{2k}\n_{L^1_tH^s_x}\leq T^2 \n u\n^{4k+2}_{L^\infty_tH^s_x},$$
so, fixing an upper bound for $T$, $T\leq 1$ for instance (but not necessarily),
$$\n v\n_{X^s_T} \leq Ce^{C' \left(\frac{T^{\alpha+1}}{\alpha+1}+T\right)}\left( \n u_0\n_s +\int_0^T \n u|u|^{2k}\n_s dt+\left(\int_0^T \n u|u|^{2k}\n^2_{H^s_x} dt\right)^{1/2}\right)$$
$$\leq Ce^{C' \left(\frac{1}{\alpha+1}+1\right)}\left( \n u_0\n_{H^s_x} +T \n u\n_{L^\infty_t H^s_x}^{2k+1}\right)$$
$$\leq C\n u_0\n_{H^s_x}+ C T \n u\n_{L^\infty_t H^s_x}^{2k+1}$$
$$\leq C\n u_0\n_{H^s_x}+ C T \n u\n_{X^s_T}^{2k+1}.$$
where, with some abuse of notations, $C$ is a new suitable constant.\\
From the previous estimate we get that $\Phi$ sends $X^s_T$ into itself. Moreover, let $R$ be $R=\frac C 2 \n u_0\n_{H^s_x}$, then, once again from the previous estimate, for all $u\in B_R\subset X_T^s$ (where $B_R$ denotes the ball of radius $R$ in $X_T^s$) we have
$$\n \Phi(u) \n_{X^s_T}=\n v\n^2_{X^s_T}\leq R/2 +CT R^{2k+1},$$
which gives, by choosing $T<1$ sufficiently small so that $CT R^{2k+1}<R/2$, that $\Phi$ sends $B_R$ into $B_R$.

What is left now is to prove that $\Phi$ is a contraction. We then consider $v:=\Phi(u)$ and $w:=\Phi(u')$, and have

$$\n v-w\n_{X^s_T}=\n \int_0^t W_\alpha(t,\tau) (u|u|^{2k}-u'|u'|^{2k})d\tau\n_{X^s_T}.$$
By the previous argument applied to $v-w$, which, in particular, is the solution of the linear problem with $f=u|u|^{2k}-u'|u'|^{2k}$ and initial datum $0$, we have (for $T<1$)
$$\n v-w\n_{X^s_T}\leq C\left( \int_0^T \n u|u|^{2k}-u'|u'|^{2k}\n_s dt+\left(\int_0^T \n u|u|^{2k}-u'|u'|^{2k}\n^2_s dt\right)^{1/2}\right).$$
By using the estimates used in Theorem \ref{teo1} and Theorem \ref{teo2}, we have
$$\n v-w\n_{X^s_T}\leq C T (\n u\n_{L^\infty_t H^s_x}^{2k} +\n u'\n^{2k}_{L^\infty_t H^s_x})\n u-v\n_{L^\infty_t H^s_x}$$
$$ \leq C T (\n u\n_{X^s_T}^{2k} +\n u'\n^{2k}_{X^s_T})\n u-u'\n_{X^s_T}.$$
Recalling that $v=\Phi(u)$ and $w=\Phi(u')$, we obtain, for any $u,u'\in B_R$,
$$\n \Phi(u)-\Phi(u')\n_{X^s_T}\leq C T R^{2k}\n u-u'\n_{X^s_T}.$$
Finally, eventually by taking $T$ smaller in such a way that $C T R^{2k}<1$, we conclude that $\Phi$ is a contraction, which gives, after application of the standard fixed point argument, the desired result.

\endproof

We now consider the case in which the nonlinearity $P(u,\bar{u},\nabla u, \nabla \bar{u})=\pm t^\beta \nabla_x u\cdot u^{2k}$, that is, $P(u,\bar{u},\nabla u, \nabla \bar{u})=\pm t^\beta \sum_{j=1}^n\partial_{x_j}u\cdot u^{2k}$.
To deal with this case we will need some lemmas that we will borrow from \cite{KPV} and that we recall in the Appendix (see Lemma \ref{lem612} and Lemma \ref{lem613}). Additionally, we will make use of the following lemma.

\begin{lemma}\label{estlem}
Let $f, g\in H^s_x(\mathbb{R}^n)$, $s>n/2$, such that $\langle x\rangle^{2N}f,\langle x\rangle^{2N}g\in H^{n/2+\varepsilon}$ for some $\varepsilon>0$ and $N\in\mathbb{N}$, then
$$\n \la x\ra ^{2N}fg\n_s\lesssim \n \la x\ra ^{2N}f\n^2_{n/2+\varepsilon}\n g\n^2_{s}+\n \la x\ra ^{2N}g\n^2_{n/2+\varepsilon}\n f\n^2_{s}.$$
\end{lemma}

\proof
Since for $|\xi-\eta|\leq |\eta|$ we have $\la\xi\ra =(1+|\xi-\eta+\eta|^2)^{1/2}\lesssim \la \eta \ra$, then
$$\n \la x\ra ^{2N}fg\n_s^2:=\int_{\mathbb{R}^n}\la \xi\ra^{2s} |(I-\triangle_\xi)^N\widehat{fg}(\xi)|^2 d\xi$$
$$=\int_{\mathbb{R}^n}\la \xi\ra^{2s} \left|(I-\triangle_\xi)^N\left(\int_{\mathbb{R}^n}\widehat{f}(\xi-\eta)\widehat{g}(\eta) d\eta \right)\right|^2 d\xi$$

$$ \int_{\mathbb{R}^n}\la \xi\ra^{2s} \left|(I-\triangle_\xi)^N\left(\int_{ |\xi-\eta|> |\eta|}\widehat{f}(\xi-\eta)\widehat{g}(\eta) d\eta+ \int_{ |\xi-\eta|\leq |\eta|}\widehat{f}(\xi-\eta)\widehat{g}(\eta) d\eta \right)\right|^2 d\xi$$
$$\leq \int_{\mathbb{R}^n}\la \xi\ra^{2s} \left|(I-\triangle_\xi)^N\left(\int_{|\xi-\eta|> |\eta|}\widehat{f}(\xi-\eta)\widehat{g}(\eta) d\eta \right)\right|^2 d\xi$$
$$+\int_{\mathbb{R}^n}\la \xi\ra^{2s} \left|(I-\triangle_\xi)^N\left( \int_{|\xi-\eta|\leq |\eta|}\widehat{f}(\xi-\eta)\widehat{g}(\eta) d\eta \right)\right|^2 d\xi$$
$$\lesssim \int_{\mathbb{R}^n}\la \xi\ra^{2s} \left|(I-\triangle_\xi)^N\left(\int_{ |\gamma|> |\xi-\gamma|}\widehat{f}(\gamma)\widehat{g}(\xi-\gamma) d\eta \right)\right|^2 d\xi$$
$$ +\int_{\mathbb{R}^n} \left( \int_{ |\xi-\eta|\leq |\eta|}\la \eta\ra^{s}|(I-\triangle_\xi)^N\widehat{f}(\xi-\eta)||\widehat{g}(\eta)| d\eta \right)^2 d\xi$$
$$\lesssim \int_{\mathbb{R}^n} \left(\int_{\mathbb{R}^n}\la \gamma\ra^{s}|\widehat{f}(\gamma)||(I-\triangle_\xi)^N\widehat{g}(\xi-\gamma)| d\eta \right)^2 d\xi$$
$$+\int_{\mathbb{R}^n} \left(\int_{ \mathbb{R}^n}\la \eta\ra^{s}|(I-\triangle_\xi)^N\widehat{f}(\xi-\eta)||\widehat{g}(\eta)| d\eta \right)^2 d\xi$$
$$ \lesssim \n  \, |(I-\triangle_\xi^N)\widehat{g}|\ast \la \xi\ra^s|\widehat{f}|\,\n_{L^2}^2+ \n  \,|(I-\triangle_\xi)^N\widehat{f}\,|\ast \la \xi\ra^s|\widehat{g}|\,\n_{L^2}^2$$
$$\underset{\text{Joung's ineq.}}{\lesssim} \n\la \xi\ra^s |\widehat{f}|\,|\n_{L^2}^2\n(I-\triangle_\xi)^N\widehat{g}\n^2_{L^1}+ 
\n(I-\triangle_\xi)^N\widehat{f}\n^2_{L^1} \n \la \xi\ra^s |\widehat{g}|\n_{L^2}^2$$

$$\lesssim\n f\n_s^2\int_{\mathbb{R}^n}\frac{1}{ \la \xi\ra ^{n/2+\varepsilon}}\la \xi\ra ^{n/2+\varepsilon}|(I-\triangle_\xi)^N\widehat{g}\,|\, d\xi+\n g\n_s^2\int_{\mathbb{R}^n}\frac{1}{ \la \xi\ra ^{n/2+\varepsilon}}\la \xi\ra ^{n/2+\varepsilon}|(I-\triangle_\xi)^N\widehat{f}|\, d\xi$$
$$\lesssim \n f\n_s^2 \n\la \xi\ra ^{n/2+\varepsilon} (I-\triangle_\xi)^N\widehat{g}\n_{L^2}^2+\n g\n_s^2 \n \la \xi\ra ^{n/2+\varepsilon}(I-\triangle_\xi)^N\widehat{f}\n_{L^2}$$
$$\lesssim  \n f\n_s^2 \n \la x\ra^{2N}g\n_{n/2+\varepsilon}+ \n g\n_s^2 \n \la x\ra^{2N}f\n_{n/2+\varepsilon},$$
which concludes the proof.
\endproof

We are now ready to finish the proof of Theorem \ref{teo2}.
\proof[Proof of Theorem \ref{teo2}]
First we assume that $\beta>\alpha$. Once more we consider the focusing case and write the solution of the IVP under consideration as
$$u(t,x)=W_\alpha(t)u_0+\int_0^tW_\alpha(t,\tau)\tau^\beta \nabla_xu\cdot u^{2k}d\tau.$$
We look for the solution given by the fixed point of the map 
$$\Phi(u):= W_\alpha(t)u_0+\int_0^tW_\alpha(t,\tau){\tau}^\beta\nabla_xu\cdot u^{2k}d\tau,$$
now defined on 
$$X_T^s:=\{u: [0,T]\times \mathbb{R}^n\rightarrow \mathbb{C}; \n u\n_{L^\infty_t H^s_x}<\infty,\,\left( \int_0^T\int_{\mathbb{R}^n}t^\alpha\lambda(|x|)|\Lambda^{s+1/2}u|^2 dx\,  dt \right)^{1/2}<\infty,\,$$
$$ \n \lambda(|x|)^{-1}u\n_{L^\infty_tH^{s-2N-3/2}_x}<\infty\},$$
where
$$\n u\n_{X_T^s}^2= \n u\n^2_{L^\infty_t H^s_x}+ \int_0^T\int_{\mathbb{R}^n}t^\alpha\lambda(|x|)|\Lambda^{s+1/2}u|^2 dx\,  dt + \n \lambda(|x|)^{-1}u\n_{L^\infty_tH^{s-2N-3/2}_x}^2.$$
We then call $v:=\Phi(u)$ the solution of the linear problem

\begin{equation*}
\left\{ \begin{array}{ll}
\partial_tv-it^\alpha\Delta_xv -ib(t,x)\cdot \nabla_xv=t^\beta\nabla_xu\cdot u^{2k}& \\
v(0,x)=u_0(x), & 
\end{array}\right .
\end{equation*}
and, as before, we make use of the linear smoothing estimates to prove that $\Phi$ is a contraction.
In the sequel, for shortness, we will often use the notations $\lambda:=\lambda(|x|)=\la x\ra ^{-2N}$, with $N\geq 1$ (i.e. $\sigma=2N$), and $\nabla:=\nabla_x:=\sum_{j=1}^n\partial_{x_j}$.
We have

$$\n v\n^2_{X^s_T}:=\n v\n^2_{L^\infty_t H^s_x} +\int_0^T\int_{\mathbb{R}^n}t^\alpha\lambda(|x|)|\Lambda^{s+1/2}v|^2 dx\,  dt + \n \lambda(|x|)^{-1}v\n^2_{L^\infty_tH^{s-2N-3/2}_x}$$
$$=I+II+III,$$
(recall $s>n+4N+3$), and we estimate the three terms separately.
By application of (iii) of Theorem \ref{thb} (we assume $T\leq 1$ and estimate the exponentials with exponent depending on time directly with a suitable uniform constant) we have
$$I+II=\n v\n^2_{L^\infty_t H^s_x} +\int_0^T\int_{\mathbb{R}^n}t^\alpha\lambda(|x|)|\Lambda^{s+1/2}v|^2 dx\,  dt $$
$$\lesssim \n u_0^2\n_{L^\infty_t H^s_x}+\int_0^T\int_{\mathbb{R}^n}t^{-\alpha}\lambda(|x|)^{-1}|\Lambda^{s-1/2}({t}^\beta \nabla u\cdot u^{2k})|^2 dx\,  dt$$
$$=\n u_0\n^2_{L^\infty_t H^s_x}+\int_0^T\int_{\mathbb{R}^n}t^{2\beta-\alpha}\lambda(|x|)^{-1}|\Lambda^{s-1/2}( \nabla u\cdot u^{2k})|^2 dx\,  dt$$
$$=\n u_0\n^2_{L^\infty_t H^s_x}+ II'.$$
Since $s-1/2\in 2\mathbb{N}$, then $\Lambda^{s-1/2}$ is a differential operator and, by Leibnitz rule, we have

$$\Lambda^{s-1/2} (\nabla u\cdot u^{2k})= (\Lambda^{s-1/2} \nabla u) u^{2k}
+\underset{s/2-1/4\leq|\gamma_1|<s-1/2,\,|\gamma_2|\leq s/2-1/4}{\underset{|\gamma_1|+|\gamma_2|\leq s-1/2}{\sum}}  
C_{\gamma_1,\gamma_2,s}(D^{\gamma_1}\nabla_xu) (D^{\gamma_2}u^{2k})$$
$$+ \underset{|\gamma_1|<s/2-1/4,\,|\gamma_2|> s/2-1/4}{\underset{|\gamma_1|+|\gamma_2|\leq s-1/2}{\sum}} 
C_{\gamma_1,\gamma_2,s}(D^{\gamma_1}\nabla_xu) (D^{\gamma_2}u^{2k}),$$
and

$$II'\leq \int_0^T\int_{\mathbb{R}^n}t^{2\beta-\alpha}\lambda(|x|)^{-1}|(\Lambda^{s-1/2} \nabla u) u^{2k}|^2 dx\,  dt$$
$$+ \underset{s/2-1/4\leq|\gamma_1|<s-1/2,\,|\gamma_2|\leq s/2-1/4}{\underset{|\gamma_1|+|\gamma_2|\leq s-1/2}{\sum}}
C_{\gamma_1,\gamma_2,s}\int_0^T\int_{\mathbb{R}^n}t^{2\beta-\alpha}\lambda(|x|)^{-1}|(D^{\gamma_1}\nabla_xu) (D^{\gamma_2}u^{2k})|^2 dx\,dt$$
$$+  \underset{|\gamma_1|<s/2-1/4,\,|\gamma_2|> s/2-1/4}{\underset{|\gamma_1|+|\gamma_2|\leq s-1/2}{\sum}} 
C_{\gamma_1,\gamma_2,s}\int_0^T\int_{\mathbb{R}^n}t^{2\beta-\alpha}\lambda(|x|)^{-1}|(D^{\gamma_1}\nabla_xu) (D^{\gamma_2}u^{2k})|^2 dx\,dt$$
$$= II'_a+II'_b+II'_c.$$

For $II'_a$ we have

$$II'_a\leq T^{2\beta-2\alpha} \int_0^T\int_{\mathbb{R}^n}t^{\alpha}\lambda(|x|)|\Lambda^{s-1/2} \nabla u|^2\cdot |\lambda(|x|)^{-1} u^{2k}|^2 dx\,  dt$$
$$\leq T^{2\beta-2\alpha} \left(\int_0^T\int_{\mathbb{R}^n}t^{\alpha}\lambda(|x|)|\Lambda^{s+1/2} u|^2dx\,  dt\right)\cdot \n \lambda^{-1} u^{2k}\n^2_{L^\infty_t L^\infty_x}$$
$$\leq T^{2\beta-2\alpha} \n u\n^2_{X^s_T} \,\n \lambda^{-1}u\n^2_{L^\infty_t H^{s-2N-3/2}_x} \n u\n^{4k-2}_{L^\infty_t H^s_x}$$
$$\leq T^{2\beta-2\alpha} \n u\n^{4k+2}_{X^s_T}.$$
For $II'_b$ we have

$$II'_b=\underset{s/2-1/4\leq |\gamma_1|<s-1/2,\,|\gamma_2|\leq s/2-1/4}{\underset{|\gamma_1|+|\gamma_2|\leq s-1/2}{\sum}} 
C_{\gamma_1,\gamma_2,s}\int_0^T\int_{\mathbb{R}^n}t^{2\beta-\alpha}\lambda(|x|)^{-1}|(D^{\gamma_1}\nabla_xu) (D^{\gamma_2}u^{2k})|^2 dx\,dt$$
$$= \underset{|\gamma_1|<s-1/2,\,|\gamma_2|\leq s/2-1/4}{\underset{|\gamma_1|+|\gamma_2|\leq s-1/2}{\sum}} 
C_{\gamma_1,\gamma_2,s}\int_0^T\int_{\mathbb{R}^n}t^{2\beta-\alpha}|(D^{\gamma_1}\nabla_xu)|^2 |\lambda(|x|)^{-1/2}D^{\gamma_2}u^{2k}|^2 dx\,dt$$
$$\leq \underset{|\gamma_1|<s-1/2,\,|\gamma_2|\leq s/2-1/4}{\underset{|\gamma_1|+|\gamma_2|\leq s-1/2}{\sum}} 
C_{\gamma_1,\gamma_2,s} \,T^{2\beta-\alpha}\left( \int_0^T\int_{\mathbb{R}^n}|(D^{\gamma_1}\nabla_xu)|^2 dx\,dt\right)\cdot \n \lambda(|x|)^{-1}D^{\gamma_2}u^{2k}\n^2_{L^\infty_t L^\infty_x}$$
$$\leq T^{2\beta-\alpha}\underset{|\gamma_1|<s-1/2,\,|\gamma_2|\leq s/2-1/4}{\underset{|\gamma_1|+|\gamma_2|\leq s-1/2}{\sum}} 
C_{\gamma_1,\gamma_2,s}  \n u\n_{L^\infty_t H^{\gamma_1+1}}^2\cdot  \n \lambda(|x|)^{-1}D^{\gamma_2}u^{2k}\n^2_{L^\infty_t H^{n/2+\varepsilon}_x}.$$
Note that,  denoting by $\Psi^k$  a pseudo-differential operator (with constant coefficients) of order $k$, by using Lemma \ref{lem612} we have

 $$\n \lambda(|x|)^{-1}D^{\gamma_2}u^{2k}\n^2_{L^\infty_t H^{n/2+\varepsilon}_x}\lesssim  \n D^{\gamma_2}\lambda(|x|)^{-1}u^{2k}\n^2_{L^\infty_t H^{n/2+\varepsilon}_x}+\sum_{j=1}^n\n \Psi^{|\gamma_2|-1}x_j\la x\ra^{2N-2}u^{2k}\n_{L^\infty_t H^{n/2+\varepsilon}_x}$$
 $$+\sum_{|\alpha+\beta|\leq 2N,\, |\alpha|>2,\, |\beta|\leq 2N-2}\n \Psi^{|\gamma_2|-|\alpha|}x^\beta u^{2k}\n_{L^\infty_t H^{n/2+\varepsilon}_x}$$
 $$\lesssim \n \lambda(|x|)^{-1}u^{2k}\n^2_{L^\infty_t H^{n/2+\varepsilon+|\gamma_2|}_x}+
 \sum_{j=1}^n\n x_j\la x\ra^{2N-2}u^{2k}\n_{L^\infty_t H^{n/2+\varepsilon+|\gamma_2|-1}_x}$$
 $$+\sum_{|\alpha+\beta|\leq 2N,\, |\alpha|>2,\, |\beta|\leq 2N-2}\n x^\beta u^{2k}\n_{L^\infty_t H^{n/2+\varepsilon+|\gamma_2|-|\alpha|}_x}$$
 
  $$\lesssim \n \lambda(|x|)^{-1}u^{2k}\n^2_{L^\infty_t H^{n/2+\varepsilon+|\gamma_2|}_x}+
 \sum_{j=1}^n\n\frac{ x_j}{\la x\ra^2}\la x\ra^{2N}u^{2k}\n_{L^\infty_t H^{n/2+\varepsilon+|\gamma_2|-1}_x}$$
 $$+\sum_{|\alpha+\beta|\leq 2N,\, |\alpha|>2,\, |\beta|\leq 2N-2}\n \frac{x^\beta}{\la x\ra^{2N}} \la x\ra^{2N} u^{2k}\n_{L^\infty_t H^{n/2+\varepsilon+|\gamma_2|-|\alpha|}_x}$$
 $$\lesssim  \n \lambda(|x|)^{-1}u^{2k}\n^2_{L^\infty_t H^{n/2+\varepsilon+|\gamma_2|}_x}$$
 $$\lesssim  \n \lambda(|x|)^{-1}u^{}\n^2_{L^\infty_t H^{n/2+\varepsilon+|\gamma_2|}_x}\n u\n^{4k-2}_{L^\infty_t H^{n/2+\varepsilon+|\gamma_2|}_x},$$
where we used the $H^s$-boundedness of $\frac{x_j}{\la x\ra^{2}}$ and $\frac{x^\beta}{\la x\ra^{2N}}$ as pseudo-differential operators of order 0, together with Sobolev inequalities.

By using the previous estimate in $II'_b$ and using  $\n \cdot\n_{H^{n/2+\varepsilon+|\gamma_2|}}\leq \n \cdot\n_{H^{s-2N-3/2}}$ for $\varepsilon$ sufficiently small such that $s/2+n/2+\varepsilon-1/4\leq s-2N-3/2$ (recall $s\geq n+4N+3$), we get

$$II'_b\lesssim T^{2\beta-\alpha} \n u\n_{X^s_T}^{4k+2}.$$

For $II'_c$, repeating the steps in the estimate of $II'_b$, we have

$$II'_c=\underset{|\gamma_1|<s/2-1/4,\,|\gamma_2|> s/2-1/4}{\underset{|\gamma_1|+|\gamma_2|\leq s-1/2}{\sum}} 
C_{\gamma_1,\gamma_2,s}\int_0^T\int_{\mathbb{R}^n}t^{2\beta-\alpha}\lambda(|x|)^{-1}|(D^{\gamma_1}\nabla_xu) (D^{\gamma_2}u^{2k})|^2 dx\,dt$$
$$ \leq \underset{|\gamma_1|<s/2-1/4,\,|\gamma_2|> s/2-1/4}{\underset{|\gamma_1|+|\gamma_2|\leq s-1/2}{\sum}} 
C_{\gamma_1,\gamma_2,s}\,T^{2\beta-\alpha}\left( \int_0^T\int_{\mathbb{R}^n}|(D^{\gamma_2}u^{2k})|^2 dx\,dt\right)\cdot 
\n \lambda(|x|)^{-1}D^{\gamma_1}\nabla u\n^2_{L^\infty_t L^\infty_x}$$
$$\leq T^{2\beta-\alpha} \underset{|\gamma_1|<s/2-1/4,\,|\gamma_2|> s/2-1/4}{\underset{|\gamma_1|+|\gamma_2|\leq s-1/2}{\sum}} C_{\gamma_1,\gamma_2,s} \n u\n^{4k}_{L^\infty_t H^{\gamma_2}_x}\n \lambda(|x|)^{-1} D^{\gamma_1+1}u\n^2_{L^\infty_t H^{n/2+\varepsilon}}$$
$$\underset{\text{by Lemma \ref{lem612}}}{\lesssim} T^{2\beta-\alpha} \n u\n^{4k+2}_{X_T^s},$$

and finally, for $T\leq 1$,

$$I+II\leq C T^{2\beta-2\alpha} \n u\n_{X_T^s}^{4k+2}.$$

To estimate $III$ we use Lemma \ref{lem613}, so we have
 
 $$III=\n \lambda(|x|)^{-1}v\n^2_{L^\infty_tH^{s-2N-3/2}_x} $$
 $$\lesssim \n \lambda(|x|)^{-1} W_\alpha(t)u_0\n^2_{L^\infty_tH^{s-2N-3/2}_x}+\n\int_0^t  \lambda(|x|)^{-1}  W_\alpha(t,\tau){\tau}^\beta \nabla u \cdot u^{2k}d\tau\n^2_{L^\infty_tH^{s-2N-3/2}_x}$$
$$\underset{\text{by Lemma \ref{lem613} and Minkowski}}{\lesssim} C(1+T^{2N})^2\n \lambda^{-1} u_0\n^2_{H^{s-3/2}_x} +$$
$$ \left(\int_0^T \n  \lambda(|x|)^{-1}  W_\alpha(t,\tau){\tau}^\beta \nabla u\cdot u^{2k}\n_{L^\infty_{\tau}H^{s-2N-3/2}_x} dt\right)^2$$
$$\lesssim  C(1+T^{2N})^2\n \lambda^{-1} u_0\n^2_{H^{s-3/2}_x} + \left(T \underset{0\leq \tau\leq t\leq T}{\sup}\n  \lambda(|x|)^{-1}  W_\alpha(t,\tau){\tau}^\beta \nabla u\cdot u^{2k}\n_{H^{s-2N-3/2}_x}\right)^2$$
$$\underset{\text{by Lemma \ref{lem613}}}{\lesssim}  C(1+T^{2N})^2\n \lambda^{-1} u_0\n^2_{H^{s-3/2}_x} + \left(CT(1+T^{2N}) \underset{0\leq  t\leq T}{\sup}\n  \lambda^{-1}  {t}^\beta \nabla u\cdot u^{2k}\n_{H^{s-3/2}_x}\right)^2$$
$$\lesssim C(1+T^{2N})^2\n \lambda^{-1} u_0\n^2_{H^{s-3/2}_x} + CT^{2(\beta+1)}(1+T^{2N})^2\n \lambda^{-1}\nabla u\cdot u^{2k}\n^2_{L^\infty_t H^{s-3/2}_x}.$$
$$\underset{\text{by Lemma \ref{estlem}}}{\lesssim}C(1+T^{2N})^2\n \lambda^{-1} u_0\n^2_{H^{s-3/2}_x} + CT^{2(\beta+1)}(1+T^{2N})^2
(\n \lambda^{-1}\nabla u\n_{L^\infty_t H^{n/2+\varepsilon}_x}^2\n u^{2k}\n^2_{L^\infty_t H^{s-3/2}_x}$$
$$+\n \lambda^{-1}u^{2k}\n_{L^\infty_t H^{n/2+\varepsilon}_x}^2\n \nabla u\n_{L^\infty_t H^{s-3/2}_x}^2).$$

By Lemma \ref{lem612} we have (recall $\nabla:=\sum_{j=1}^n\partial_{x_j}$)
$$\lambda^{-1}\nabla (u)=\nabla (\lambda^{-1}u)+\sum_{j=1}^n D_{x_j}(x_j\la x\ra^{2N-2} u),$$
therefore, 
$$\n \lambda^{-1}\nabla u\n_{L^\infty_t H^{n/2+\varepsilon}_x}\lesssim
\n \nabla \lambda^{-1}u\n_{L^\infty_t H^{n/2+\varepsilon}_x}+ \n\sum_{j=1}^nD_{x_j}x_j\la x\ra^{2N-2} u\n_{L^\infty_t H^{n/2+\varepsilon}_x}$$
$$=\n \nabla \lambda^{-1}u\n_{L^\infty_t H^{n/2+\varepsilon}_x}+ \n\sum_{j=1}^nD_{x_j}\frac{x_j}{\la x\ra^2}\la x\ra^{2N} u\n_{L^\infty_t H^{n/2+\varepsilon}_x}$$
$$\lesssim \n \lambda^{-1}u\n_{L^\infty_t H^{n/2+\varepsilon+1}_x}\lesssim \n \lambda^{-1}u\n_{L^\infty_t H^{n/2+\varepsilon+1}_x},$$
since $\frac{x_j}{\la x\ra^{2}}\in S^0$ and, for $\varepsilon$ sufficiently small, $n/2+\varepsilon+1\leq s-2N-3/2$.
We then have

$$III\lesssim C(1+T^{2N})^2\n \lambda^{-1} u_0\n^2_{H^{s-3/2}_x} + CT^{2(\beta+1)}(1+T^{2N})^2\times$$
$$\times \left(\n\lambda^{-1}u\n_{L^\infty_t H^{n/2+\varepsilon+1}_x}^2\n u\n^{4k}_{L^\infty_tH^s_x}+\n\lambda^{-1}u\n^2_{L^\infty_tH^{n/2+\varepsilon}}\n u\n^{4k-2}_{L^\infty_tH^{s}_x } \n u\n_{L^\infty_t H^s_x}^2 \right)$$
$$\leq C(1+T^{2N})^2\n \lambda^{-1} u_0\n^2_{H^{s-3/2}_x}+ CT^{2(\beta+1)}(1+T^{2N})^2\n u\n_{X^s_T}^{4k+2}$$

From the previous estimates we get
$$\n v\n^2_{X^s_T}=I+II+III$$
$$\lesssim \n u_0\n^2_{L^\infty_t H^s_x}+ C T^{2\beta-2\alpha}\n u\n^{4k+2}_{X^s_T}+C(1+T^{2N})^2\n \lambda^{-1} u_0\n^2_{H^{s-3/2}_x} + CT^{2(\beta+1)}(1+T^{2N})^2\n u\n_{X^s_T}^{4k+2}.$$
$$\underset{T\leq 1}{\leq} C(\n u_0\n^2_{L^\infty_t H^s_x}+ \n \lambda^{-1} u_0\n^2_{H^{s-3/2}_x} )+ C(T^{2\beta-2\alpha}+T^{2(\beta+1)}(1+T^{2N})^2)\n u\n_{X^s_T}^{4k+2},$$
where, recall, $\beta>\alpha$.
Finally we have, with new suitable constants that we keep denoting simply $C$,
$$\n\Phi(u)\n_{X^s_T}=\n v\n_{X^s_T}\leq C(\n u_0\n_{L^\infty_t H^s_x}+ \n \lambda^{-1} u_0\n_{H^{s-3/2}_x} )+ C(T^{2\beta-2\alpha}+T^{2(\beta+1)}(1+T^{2N})^2)^{1/2}\n u\n_{X^s_T}^{2k+1},$$
hence, by choosing $R= C/2(\n u_0\n_{L^\infty_t H^s_x}+ \n \lambda^{-1} u_0\n_{H^{s-3/2}_x} )$ and $T$ sufficiently small such that $C(T^{2\beta-2\alpha}+T^{2(\beta+1)}(1+T^{2N})^2)^{1/2} R^{2k}<1/2$, we get that $\Phi$ sends the ball $B_R\subset X_T^s$ into itself. 

What is left to prove to conclude the proof is that $\Phi$ is a contraction.
We then consider, as in the proof of Theorem \ref{teo1}, $v:=\Phi(u)$ and $w:=\Phi(u')$ as the solutions of the linear IVP with the same initial datum $u_0$. By application of the linear smoothing estimates on $v-w$ we have
$$\n v-w\n^2_{X_T^s}\lesssim \n v-w\n^2_{L^\infty_t H^s_x} +\int_0^T\int_{\mathbb{R}^n}t^\alpha\lambda(|x|)|\Lambda^{s+1/2}(v-w)|^2 dx\,  dt + \n \lambda(|x|)^{-1}(v-w)\n^2_{L^\infty_tH^{s-2N-3/2}_x}$$
$$\lesssim \int_0^T\int_{\mathbb{R}^n}t^{2\beta-\alpha}\lambda^{-1}(|x|)|\Lambda^{s-1/2}(\nabla u\cdot u^{2k}-\nabla u'\cdot {u'}^{2k})|^2 dx\,  dt + \n \lambda(|x|)^{-1}(v-w)\n^2_{L^\infty_tH^{s-2N-3/2}_x}.$$
$$\lesssim \int_0^T\int_{\mathbb{R}^n}t^{2\beta-\alpha}\lambda^{-1}(|x|)|\Lambda^{s-1/2}\nabla (u^{2k+1}- {u'}^{2k+1})|^2 dx\,  dt + \n \lambda(|x|)^{-1}(v-w)\n^2_{L^\infty_tH^{s-2N-3/2}_x}$$
$$= IV+V.$$
For the term $IV$ we proceed like in the estimate of $II'_a$, $II'_b$ and $II'_c$ (recall $s-1/2\in 2\mathbb{N}$, so $\Lambda^{s-1/2}$ is a differential operator on which Leibnitz rule applies) and have

$$IV=\int_0^T\int_{\mathbb{R}^n}t^{2\beta-\alpha}\lambda^{-1}(|x|)|\Lambda^{s-1/2}\nabla (u^{2k+1}- {u'}^{2k+1})|^2 dx\,  dt$$
$$\leq T^{2\beta-2\alpha} \int_0^T\int_{\mathbb{R}^n}t^{\alpha}\lambda^{-1}(|x|)|\nabla \Lambda^{s-1/2}(u^{2k+1}-{u'}^{2k+1})|^2 dx\,  dt$$
$$\leq T^{2\beta-2\alpha} \int_0^T\int_{\mathbb{R}^n}t^{\alpha}\lambda^{-1}(|x|)\left|\nabla\Lambda^{s-1/2}\left((u-u')\sum_{j=0}^{2k} u^{2k-j}{u'}^j\right)\right|^2 dx\,  dt$$
$$\lesssim T^{2\beta-2\alpha} \n u-u'\n^2_{X^s_T}\n \lambda^{-1}\sum_{j=0}^{2k} u^{2k-j}{u'}^j\n^2_{L^\infty_tL^\infty_x}$$
$$ +T^{2\beta-\alpha}\n u-u'\n_{X^s_T}^2\n \lambda^{-1}\sum_{j=0}^{2k} u^{2k-j}{u'}^j\n^2_{L^\infty_tH^{s-2N-3/2}_x}$$
$$ +T^{2\beta-\alpha}\n \lambda^{-1}(u-u')\n_{L^\infty_tH^{s-2N-3/2}_x}^2\n \sum_{j=0}^{2k} u^{2k-j}{u'}^j\n^2_{L^\infty_t H^s_x}$$

$$\underset{T\leq 1}{\lesssim} T^{2\beta-2\alpha} \n u-u'\n^2_{X^s_T}\sum_{j=0}^{2k}\n \lambda^{-1} u^{2k-j}{u'}^j\n^2_{L^\infty_tH^{s-2N-3/2}_x}$$
$$+T^{2\beta-\alpha}\n u-u'\n_{X^s_T}^2 \sum_{j=0}^{2k} \n u^{2k-j}\n^2_{L^\infty_t H^s_x} \n{u'}^j\n^2_{L^\infty_t H^s_x}$$

$$\lesssim T^{2\beta-2\alpha} \n u-u'\n^2_{X^s_T}\sum_{j=0}^{2k}\n \lambda^{-1} u^{2k-j}\n_{L^\infty_tH^{s-2N-3/2}_x}\n \lambda^{-1} {u'}^{j}\n^2_{L^\infty_tH^{s-2N-3/2}_x}$$
$$+T^{2\beta-\alpha}\n u-u'\n_{X^s_T}^2(\n u\n_{X_T^s}+\n u'\n_{X_T^s})^{4k}$$
$$\lesssim T^{2\beta-2\alpha} \n u-u'\n^2_{X^s_T} \left( \n \lambda^{-1}u^{2k}\n^2_{L^\infty_tH^{s-2N-3/2}_x}
+\sum_{j=1}^{2k}\n  u^{2k-j}(\lambda^{-1}{u'}^j)\n^2_{L^\infty_tH^{s-2N-3/2}_x}\right)$$
$$+T^{2\beta-\alpha}\n u-u'\n_{X^s_T}^2(\n u\n_{X_T^s}+\n u'\n_{X_T^s})^{4k}$$

$$\lesssim T^{2\beta-2\alpha} \n u-u'\n^2_{X^s_T} \left(\n u\n^{4k}_{X^s_T}+\sum_{j=1}^{2k}\n u\n^{4k-2j}_{X^s_T}\n u'\n^{2j}_{X^s_T}\right)+T^{2\beta-\alpha}\n u-u'\n_{X^s_T}^2(\n u\n_{X_T^s}+\n u'\n_{X_T^s})^{4k}$$
$$\lesssim T^{2\beta-2\alpha} \n u-u'\n^2_{X^s_T} (\n u\n_{X^s_T}+\n u'\n_{X^s_T})^{4k}.$$
where we estimated the sum $\sum_{j=1}^{2k}\n \lambda^{-1} u^{2k-j}{u'}^j\n^2_{L^\infty_tH^{s-2N-3/2}_x}$ by decoupling each term in the form $u^{2k-j}(\lambda^{-1}u') {u'}^{j-1}$ and then by using Sobolev inequalities.

For the term $V$ we use the procedure used in the estimate of $III$ above (once again, we make use of Lemma \ref{lem612} and Lemma \ref{lem613}) and some strategies used for $IV$, and finally we obtain
$$\n v-w\n_{X_T^s}\leq  C(T^{2\beta-2\alpha}+T^{2(\beta+1)}(1+T^{2N})^2)^{1/2}(\n u\n_{X^s_T}+\n u'\n_{X^s_T})^{2k}\n u-u'\n_{X_T^s}.$$
By taking $u,u'\in B_R$, and eventually by taking the time $T$ smaller, we can conclude that
$$\n \Phi(u)-\Phi(u')\n_{X_T^s}= \n v-w\n_{X_T^s}\leq C \n u-u'\n_{X_T^s},$$
with $C<1$, so $\Phi$ is a contraction. After application of the standard fixed point argument the result follows.

\medskip
We now assume that $\beta=\alpha$. In this case we are not able to use the time factor $T^{\beta-\alpha}$ in order to obtain a contraction unless we assume that the initial data is small. This is a similar issue to the one faces in \cite{KPV1} and  \cite{KPRV}, particularly in the  part  involving the norm $\lambda_3$ in page 479. It was resolved by using a version of a mean value theorem in the time variable. More precisely  we can write
$$t^\alpha |u|^{2k}\nabla_xu= t^\alpha (|u|^{2k}-|u_0|^{2k})\nabla_xu+t^\alpha |u_0|^{2k}\nabla_xu=\nabla_xut^\alpha\int_0^t(|u|^{2k})'(s)ds+t^\alpha|u_0|^{2k}\nabla_xu$$
Then in the course of the proof $\int_0^t(|u|^{2k})'(s)ds$ will bring down $t$ which together with $t^\alpha$ will give $t^\beta, \, \beta=\alpha+1$, while $t^\alpha|u_0|^{2k}\nabla_xu$ will be incorporated in the $b(x,t)$ term. We do not write down all the details since there are essentially contained in the references mentioned earlier.

\endproof

\section{Final remarks}
We conclude this paper with some remarks about the general case $b\not\equiv 0$. However we recall, once again, that the results proved in what we call \textit{general case} are still true for the particular case $b\equiv0$.\\
\vspace{0,2cm}

\noindent 1. We assumed, in the case $b\not\equiv 0$, that for all $j=1,...,n$, $b_j$ is such that  there exist $\sigma>1$ for which
$$ |\mathsf{Im}\,\partial_x^\gamma b_j(t,x)|, |\mathsf{Re}\,\partial_x^\gamma b_j(t,x)|\lesssim t^\alpha \la x\ra^{-\sigma-|\gamma|}.$$ 
The first condition on the real part of $b_j$ is natural (see also \cite{CR,M}) since it is needed to be able to apply the Sharp G\aa rding inequality which is the key point even in the proof of the local well-posedness of the linear problem (since it gives the control on the first order term $b\cdot\nabla_x$). 
Instead, the second condition, imposed both on the real and on the imaginary part of $b_j$, is needed in order to have the control $\n r_{s-1}(t,x,D)u\n_0\leq Ct^\alpha \n u\n_{s-1}$ for the error term $r_{s-1}$ in \eqref{reminder} and get the smoothing estimates needed to deal with the nonlinear problem. 
Finally the condition on $\sigma$, that is $\sigma>1$, which is imposed in the nondegenerate case as well, is required in order to avoid a loss of derivatives of the solution (see \cite{D}).

\vspace{0,1cm}
\noindent 2. The nonlinear term $\nabla u\cdot u^{2k}$ is chosen for convenience but it can be generalized. For instance, like in \cite{KPRV}, one can consider nonlinearities given by polynomials in $u$ and $\nabla u$ and their complex conjugates.\\

\vspace{0,1cm}
\noindent 3. Possibly by using the techniques  in \cite{KPRV, M, D}, one can obtain the same smoothing and local well-posedness results for the equation
$$i\partial_tu+t^\alpha \Delta_xu+b(t,x)\cdot\nabla_x u+c(t,x)\cdot\nabla_x \bar{u}$$
assuming on the term $c(t,x)$ suitable conditions (possibly similar to that assumed for $b(t,x)$). In this case the equation should be reduced to a systems, that, after diagonalization, satisfies the desired smoothing properties from which the local well-posedness follows.\\
\vspace{0,1cm}

\noindent 4. The results proved in this paper are likely still valid for equations of the form
$$i\partial_tu+g(t) \Delta_x+b(t,x)\cdot\nabla_x,$$
provided that $g$ satisfies suitable properties, as, for example, $g$  having constant sign, vanishing at $t=0$ and such that $|\partial^\gamma_x b(t,x)|\lesssim |g(t)|\lambda(|x|)$.

\appendix
\section{}
In this section we shall recall the statement of some key lemmas we used throughout the paper.

Before giving the statement of the first lemma, which is also the crucial one, that is Doi's lemma (Lemma 2.3 in \cite{D}), we state below the conditions needed to apply this result.

According with the notation used by Doi in \cite{D}, we shall denote by (B1), (B2) and (A6) the following conditions:\\

\textit{Let $a^w(t,x,\xi)$ be the Weyl symbol of a pseudo-differential operator $A=A(t,x,D_x)$(see \cite{Ho3}). We shall say that $a^w:=a$ satisfies $(B1)$, $(B2)$ and $(A6)$ if}\\

\begin{itemize}
 \item [\textbf{(B1)}] $a(t,x,\xi)=ia_2(x,\xi)+a_1(t,x,\xi)+a_0(t,x,\xi)$, \textit{where} $a_2\in S^2_{1,0}$\textit{ is real-valued and} $a_j\in S^j_{1,0}$, for $j=0,1$;\\
   
  \item [\textbf{(B2)}] $|a_2(x,\xi)|\geq \delta |\xi|^2$ \textit{with} $x\in\mathbb{R}^n$, $|\xi|^2\geq C$, \textit{and} $\delta, C>0$;\\
  
  \item [\textbf{(A6)}] \textit{There exists a real-valued function} $q\in C^\infty (\mathbb{R}^n\times\mathbb{R}^n)$ \textit{such that, with }$C_{\alpha\beta}, C_1, C_2>0,$

\begin{equation*}
|\partial^\alpha_\xi\partial_x^\beta q(x,\xi)|\leq C_{\alpha \beta}\langle x\rangle \langle \xi\rangle^{-|\alpha|},\quad  x,\xi\in\mathbb{R}^n,
\end{equation*}

\begin{equation*}
   H_{a_2}q(x,\xi)=\{a_2,q\}(x,\xi) \geq C_1|\xi|-C_2,\quad x,\xi\in\mathbb{R}^n,
\end{equation*}
\end{itemize}
\textit{where we denoted by $S^j_{1,0}=S^j_{\rho=1,\delta=0}=:S^j$ the standard class of pseudo-differential symbols of order $j$, and by $\{\cdot,\cdot\}$  the Poisson bracket.}

\begin{lemma}[Doi \cite{D}, Lemma 2.3] \label{Doi}
Assume $(B1)$, $(B2)$ and $(A6)$. Let $\lambda(s)$ be a positive non increasing function in $C([0,\infty))$. Then\\

\begin{enumerate}
  \item If $\lambda\in L^1([0,\infty))$ there exists a real-valued symbol $p\in S^0$ and $C>0$ such that
\begin{equation}
\label{Doi1}
  H_{a_2}p\geq \lambda(|x|)|\xi|-C,\quad x,\xi\in\mathbb{R}^n;
\end{equation}
  \item If $\int_0^t \lambda(\tau) d\tau\leq C \log (t+1)+C'$, $t\geq 0$, $C,C'>0$, then there exists a real-valued symbol $p\in S^0_1(\log \langle \xi\rangle)$ such that
\begin{equation}
\label{Doi2}
  H_{a_2}p\geq \lambda(|x|)|\xi|-C_1\log \langle \xi\rangle-C_2,\quad x,\xi\in\mathbb{R}^n.
\end{equation}
\end{enumerate}
\end{lemma}

\begin{remark}\label{remDoi}
We remark that, by taking $\lambda'(|x|)=C'\lambda(|x|)$ in Doi's lemma, where $C'$ is any positive constant and $\lambda$ is as in Lemma \ref{Doi}, then we get that there exists a real-valued symbol $p\in S^0$ and a constant $C>0$ such that
\begin{equation}
\label{Doi1.1}
H_{a_2}p\geq C'\lambda(|x|)|\xi|-C,\quad x,\xi\in\mathbb{R}^n.
\end{equation}
\end{remark}

We conclude the section by giving other two useful lemmas taken from \cite{KPRV}.

\begin{lemma}[Lemma 6.1.2 of \cite{KPRV}]
\label{lem612}
Let $p\in S^m_{0,1}$, $N\in\mathbb{N}$. Then
$$(1+|x|^2)^N\Psi_p f=\Psi_p\left[ (1+|x|^2)^N f \right]+2N \sum_j \Psi_{i\partial_{\xi_j}p}\left[ x_j(1+|x|^2)^{N-1} f \right] $$
$$+\underset{|\alpha+\beta|\leq N,\,|\alpha|\geq 2,\, |\beta|\leq 2N-2}{\sum}c_{\alpha\beta}\Psi_{\partial^\alpha_{\xi} p}\left[ x^\beta f\right],$$
where $\Psi_a$ stands for the pseudo-differential operator with symbol $a$.
\end{lemma}

\begin{lemma}[Lemma 6.1.3 of \cite{KPRV}]
\label{lem613}
Let $N\in\mathbb{N}$ and $s\in \mathbb{R}$. Suppose $\la x\ra^{2N}u_0\in H^{s+2N}$. Then
$$\sup_{0\leq t\leq T}\n \la x \ra^{2N}W_1(t)u_0\n_{H^s}^2\leq \sum_{j=0}^{2N}c_jT^{j}\n \la x\ra^{2N-j}u_0\n^2_{H^{s+j}}$$
and
$$\sup_{0\leq t\leq T}\n \la x \ra^{2N}W_1(t)u_0\n_{H^s}^2\leq c(1+T^{2N})\n \la x\ra^{2N}u_0\n^2_{H^{s+2N}},$$

where $W_1$ denotes the solution operator of (6.1) in \cite{KPRV} with $f=0$.
\end{lemma}
We remark that Lemma \ref{lem613} still works in our case where the operator $W_1$ will be the solution operator of the homogeneous IVP associated with $\LL_\alpha$ and that we denoted by $W_\alpha(t):=W_\alpha(t,0)$ (for details see the proof in \cite{KPRV} pag.474).

\section*{Acknowledgement}
We wish to thank Alberto Parmeggiani for useful discussions and suggestions which helped to improve the present paper.

\end{document}